\documentclass{article}
\usepackage{amsmath,amssymb,amsthm,eucal}
\usepackage{graphics}

\newtheorem{theorem}{Theorem}
\newtheorem{lemma}{Lemma}
\newtheorem{proposition}{Proposition}
\newtheorem{definition}{Definition}
\newtheorem{corollary}{Corollary}
\newtheorem{example}{Example}
\newtheorem{remark}{Remark}

\numberwithin{theorem}{section}
\numberwithin{definition}{section}
\numberwithin{lemma}{section}
\numberwithin{corollary}{section}
\numberwithin{equation}{section}
\numberwithin{proposition}{section}
\numberwithin{example}{section}
\numberwithin{remark}{section}
\numberwithin{figure}{section}

\def\es{\varnothing}
\def\SB{\subseteq}

\def\aa{\alpha}
\def\bb{\beta}
\def\gg{\gamma}

\def\ll{\lambda}

\def\bm{\boldsymbol m}
\def\bn{\boldsymbol n}
\def\bp{\boldsymbol p}

\def\eq{\Leftrightarrow}

\def\POW{\mathfrak P}
\def\POWF{\mathfrak P_{\text{\sc f}}}

\def\FFF{{\cal F}}
\def\CCC{{\cal C}}
\def\GGG{{\cal G}}
\def\HHH{{\cal H}}
\def\LLL{{\cal L}}
\def\OOO{{\cal O}}
\def\SSS{{\cal S}}
\def\TTT{{\cal T}}

\def\QQQ{{\cal Q}}

\def\Zee{\mathbb{Z}}

\begin{document}
\title{Media Theory: Representations and Examples}
	\author{Sergei~Ovchinnikov \\
	Mathematics Department\\
	San Francisco State University\\
	San Francisco, CA 94132\\
	sergei@sfsu.edu} 
\date{\today} 
\maketitle

\begin{abstract}
In this paper we develop a representational approach to media theory. We construct representations of media by well graded families of sets and partial cubes and establish the uniqueness of these representations. Two particular examples of media are also described in detail.

\end{abstract}

\section{Introduction}

A medium is a semigroup of transformations on a possibly infinite set of states, constrained by four axioms which are recalled in Section~2 of this paper. This concept was originally introduced by Jean--Claude Falmagne in his 1997 paper~\cite{jF97} as a model for the evolution of preferences of individuals (in a voting context, for example). As such it was applied to the analysis of opinion polls~\cite{mRjFbG99} (for closely related papers, see~\cite{jF96,jFjD97,jFmRbG97}). As shown by Falmagne and Ovchinnikov~\cite{jFsO02} and Doignon and Falmagne~\cite{jDjF97}, the concept of a medium provides an algebraic formulation for a variety of geometrical and combinatoric objects. The main theoretical developments so far can be found in~\cite{jF97,jFsO02,sOaD00} (see also Eppstein and Falmagne's paper~\cite{dE04} in this volume).

The purpose of this paper is to further develop our understanding of media. We focus in particular on representations of media by means of well graded families of sets and graphs. Our approach utilizes natural distance and betweenness structures of media, graphs, and families of sets. The main results of the paper show that, in some precise sense, any medium can be uniquely represented by a well graded family of sets or a partial cube. Two examples of infinite media are explored in detail in Sections~8 and~9.

\section{Preliminaries}
In this section we recall some definitions and theorems from~\cite{jF97}.

Let $\SSS$ be a set of \emph{states}. A \emph{token} (of information) is a function $\tau:\SSS\rightarrow\SSS$. We shall use the abbreviations $S\tau=\tau(S),$ and 
$S\tau_1\cdots\tau_n=\tau_n[\ldots[\tau_1(S)]]$ for the function composition. We denote by $\tau_0$ the identity function on $\SSS$ and suppose that $\tau_0$ is not a token. Let $\TTT$ be a set of tokens on $\SSS$. The pair $(\SSS,\TTT)$ is called a \emph{token system}. Two distinct states $S,T\in\SSS$ are \emph{adjacent} if $S\tau =T$ for some token $\tau\in\TTT$. To avoid trivialities, we assume that $|\SSS|>1$.

A token $\tau'$ is a \emph{reverse} of a token $\tau$ if for all distinct
$S,V\in\SSS$ 
$$
S\tau=V\quad\eq\quad V\tau'=S.
$$
A finite composition $\bm=\tau_1\cdots\tau_n$ of not necessarily distinct tokens $\tau_1,\ldots,\tau_n$ such that $S\bm=V$ is called a \emph{message producing} $V$ \emph{from} $S$. We write $\ell(\bm)=n$ to denote the \emph{length} of $\bm$. The \emph{content} of a message
$\bm=\tau_1\ldots\tau_n$ is the set $\CCC(\bm)$ of its distinct tokens.  Thus, $|\CCC(\bm)|\leq\ell(\bm)$. A message $\bm$ is \emph{effective} (resp. \emph{ineffective}) for a state $S$ if $S\bm\neq S$ (resp. $S\bm=S$). A message $\bm=\tau_1\ldots\tau_n$ is \emph{stepwise effective} for $S$ if 
$$
S\tau_1\ldots\tau_k\neq S\tau_0\ldots\tau_{k-1},\qquad 1\leq k\leq n.
$$
A message is called \emph{consistent} if it does not contain both a token and its reverse, and \emph{inconsistent} otherwise. A message which is both consistent and stepwise effective for some state $S$ is said to be \emph{straight} for $S$. A message $\bm=\tau_1\ldots\tau_n$ is \emph{vacuous} if the set of indices $\{1,\ldots,n\}$ can be partitioned into pairs $\{i,j\},$ such that one of $\tau_i,\tau_j$ is a reverse of the other. Two messages $\bm$ and $\bn$ are \emph{jointly consistent} if $\boldsymbol{mn}$ (or, equivalently, $\boldsymbol{nm}$) is consistent. 

The next definition introduces the main concept of media theory.

\begin{definition} \label{D:medium}
A token system is called a \emph{medium} if the following axioms are satisfied.
\begin{enumerate}
\item[]
\begin{enumerate}
\item[{\rm [M1]}] Every token $\tau$ has a unique reverse, which we denote by $\tilde{\tau}$.
\item[{\rm [M2]}] For any two distinct states $S,V,$ there is a consistent message transforming $S$ into $V$.
\item[{\rm [M3]}] A message which is stepwise effective for some state is ineffective for that state if and only if it is vacuous.
\item[{\rm [M4]}] Two straight messages producing the same state are jointly consistent.
\end{enumerate}
\end{enumerate}
\end{definition}

It is easy to verify that \lbrack M2\rbrack\ is equivalent to the following
axiom (cf.~\cite{jF97}, Theorem 1.7).

\begin{enumerate}
\item[]
\begin{enumerate}
\item[{\rm [M2*]}] For any two distinct states $S,V,$ there is a
straight message transforming $S$ into $V$.
\end{enumerate}
\end{enumerate}
We shall use this form of axiom [M2] in the paper.

Various properties of media have been established in~\cite{jF97}. First, we recall the concept of `content'.

\begin{definition}
Let $(\SSS,\TTT)$ be a medium. For any state $S$, the \emph{content} of $S$ is the set $\widehat{S}$ of all tokens each of which is contained in at least one straight message producing $S$. The family
$\widehat{\SSS}=\{\widehat{S}\,|\,S\in\SSS\}$ is called the \emph{content family} of $\SSS$.
\end{definition}

The following theorems present results of theorems 1.14, 1.16, and 1.17 in~\cite{jF97}. For reader's convinience, we prove these results below.

\begin{theorem} \label{Theorem1.17-1}
For any token $\tau$ and any state $S$, we have either $\tau\in\widehat{S}$ or $\tilde{\tau}\in\widehat{S}$. Consequently, $|\widehat{S}|=|\widehat{V}|$ for any two states $S$ and $V$. {\rm (}$|A|$ stands for the cardinality of the set $A$.{\rm )}
\end{theorem}

\begin{proof} 
Since $\tau$ is a token, there are two states $V$ and $W$ such that $W=V\tau$. By Axiom [M2*], there are straight messages $\bm$ and $\bn$ such that $S=V\bm$ and $S=W\bn$. By Axiom [M3], the message $\tau\bn\widetilde{\bm}$ is vacuous. Therefore, $\tilde{\tau}\in\CCC(\bn)$ or $\tilde{\tau}\in\CCC(\widetilde{\bm})$. It follows that $\tilde{\tau}\in\widehat{S}$ or $\tau\in\widehat{S}$. By Axiom [M4], we cannot have both $\tilde{\tau}\in\widehat{S}$ and $\tau\in\widehat{S}$.
\end{proof}

\begin{theorem} \label{Theorem1.16}
If $S$ and $V$ are two distinct states, with $S\bm=V$ for some straight message $\bm$, then \mbox{$\widehat{V}\setminus\widehat{S}=\CCC(\bm)$}.
\end{theorem}

\begin{proof} 
If $\tau\in\CCC(\bm)$, then $\tilde{\tau}\in\CCC(\widetilde{\bm})$. Thus, $\tau\in\widehat{V}$ and $\tilde{\tau}\in\widehat{S}$. By Theorem~\ref{Theorem1.17-1}, the latter inclusion implies $\tau\notin\widehat{S}$. It follows that $\CCC(\bm)\SB\widehat{V}\setminus\widehat{S}$.

If $\tau\in\widehat{V}\setminus\widehat{S}$, then there is a state $W$ and a straight message $\bn$ such that $V=W\bn$ and $\tau\in\CCC(\bn)$. By Axiom [M2*], there is a straight message $\bp$ such that $W=S\bp$. By Axiom [M3], the message $\bp\bn\widetilde{\bm}$ is vacuous, so $\tilde{\tau}\in\CCC(\bp\bn\widetilde{\bm})$, that is, $\tau\in\CCC(\bm\widetilde{\bn}\widetilde{\bp})$. But $\tau\notin\CCC(\widetilde{\bp})$, since $S=W\widetilde{\bp}$ and $\tau\notin\widehat{S}$, and $\tau\notin\CCC(\widetilde{\bn})$, since $\tau\in\CCC(\bn)$. Hence, $\tau\in\CCC(\bm)$, that is, $\widehat{V}\setminus\widehat{S}\SB\CCC(\bm)$.
\end{proof}

\begin{theorem} \label{Theorem1.17-2}
For any token $\tau$ and any state $S$, we have either $\tau\in\widehat{S}$ or $\tilde{\tau}\in\widehat{S}$. Moreover, 
$$
S=V\quad\eq\quad\widehat{S}=\widehat{V}.
$$
\end{theorem}

\begin{proof}
Let $\widehat{S}=\widehat{V}$ and let $\bm$ be a straigt message producing $V$ from $S$. By Theorem~\ref{Theorem1.16},
$$
\es=\widehat{V}\setminus\widehat{S}=\CCC(\bm).
$$
Thus, $S=V$.
\end{proof}

\begin{theorem} \label{Theorem1.14}
Let $\bm$ and $\bn$ be two distinct straight messages transforming some state $S$. Then $S\bm=S\bn$ if and only if $\CCC(\bm)=\CCC(\bn)$.
\end{theorem}

\begin{proof}
Suppose that $V=S\bm=S\bn$. By Theorem~\ref{Theorem1.16},
$$
\CCC(\bm)=\widehat{V}\setminus\widehat{S}=\CCC(\bn).
$$

Suppose that $\CCC(\bm)=\CCC(\bn)$ and let $V=S\bm$ and $W=S\bn$. By Theorem~\ref{Theorem1.16},
$$
\widehat{V}\Delta\widehat{S}=\CCC(\bm)\cup\CCC(\widetilde{\bm})=\CCC(\bn)\cup\CCC(\widetilde{\bn})=\widehat{W}\Delta\widehat{S},
$$
which implies $\widehat{V}=\widehat{W}$. By Theorem~\ref{Theorem1.17-2}, $V=W$.
\end{proof}

One particular class of media plays an important role in our constructions (cf.~\cite{jFsO02}).

\begin{definition} \label{D:complete-media}
A medium $(\SSS,\TTT)$ is called \emph{complete} if for any state $S\in\SSS$ and token $\tau\in\TTT$, either $\tau$ or $\tilde{\tau}$ is effective on $S$.
\end{definition}
An important example of a complete medium is found below in Theorem~\ref{CompleteMedium}. 

\section{Isomorphisms, embeddings, and token subsystems}

The purpose of combinatorial media theory is to find and examine those properties of media that do not depend on a particular structure of individual states and tokens. For this purpose we introduce the concepts of embedding and isomorphism for token systems.

\begin{definition} \label{D:embedding}
Let $(\SSS,\TTT)$ and $(\SSS',\TTT')$ be two token systems. A pair $(\aa,\bb)$ of one--to--one functions $\aa:\SSS\rightarrow\SSS'$ and $\bb:\TTT\rightarrow\TTT'$ such that 
$$
S\tau=T \quad \eq \quad \aa\left(S\right)\bb\left(\tau\right)=\aa\left(T\right)
$$
for all $S,T\in\SSS$, $\tau\in\TTT$ is called an \emph{embedding} of the token system $(\SSS,\TTT)$ into the token system $(\SSS',\TTT')$.

Token systems $(\SSS,\TTT)$ and $(\SSS',\TTT')$ are \emph{isomorphic} if there is an embedding $(\aa,\bb)$ from $(\SSS,\TTT)$ into $(\SSS',\TTT')$ such that both $\aa$ and $\bb$ are bijections.
\end{definition}

Clearly, if one of two isomorphic token systems is a medium, then the other one is also a medium.

A general remark is in order. If a token system $(\SSS,\TTT)$ is a medium and $S\tau_1=S\tau_2\not=S$ for some state $S$, then, by axiom [M3], $\tau_1=\tau_2$. In particular, if $(\aa,\bb)$ is an embedding of a medium into a medium, then $\bb(\tilde{\tau})=\widetilde{\bb(\tau)}$. We extend $\bb$ to the semigroup of messages  by defining $\bb(\tau_1\cdots\tau_k)=\bb(\tau_1)\cdots\bb(\tau_k)$. Clearly, the image $\bb(\bm)$ of a straight message $\bm$ is a straight message.

Let $(\SSS,\TTT)$ be a token system and $\QQQ$ be a subset of $\SSS$ consisting of more than two elements. The restriction of a token $\tau\in\TTT$ to $\QQQ$ is not necessarily a token on $\QQQ$. In order to construct a medium with the set of states $\QQQ$, we introduce the following concept.

\begin{definition} \label{D:restriction}
Let $(\SSS,\TTT)$ be a token system, $\QQQ$ be a nonempty subset of $\SSS$, and
$\tau\in\TTT$. We define a \emph{reduction} of $\tau$ to $\QQQ$ by 
$$
S\tau_{\QQQ} = \begin{cases}
	S\tau & \text{if $S\tau\in\QQQ$,} \\
	S & \text{if $S\tau\notin\QQQ$,}
\end{cases}
$$
for $S\in\QQQ$. A token system $(\QQQ,\TTT_{\QQQ})$ where
$\TTT_{\QQQ}=\{\tau_{\QQQ}\}_{\tau\in\TTT} \setminus\{\tau_0\}$ is the set of all distinct reductions of tokens in $\TTT$ to $\QQQ$ different from the identity function $\tau_0$ on $\QQQ$, is said to be the \emph{reduction} of $(\SSS,\TTT)$ to $\QQQ$.

We call $(\QQQ,\TTT_{\QQQ})$ a \emph{token subsystem} of $(\SSS,\TTT)$. If both $(\SSS,\TTT)$ and $(\QQQ,\TTT_{\QQQ})$ are media, we call
$(\QQQ,\TTT_{\QQQ})$ a \emph{submedium} of $(\SSS,\TTT)$.
\end{definition}

A reduction of a medium is not necessarily a submedium of a given medium. Consider, for instance, the following medium:

{\begin{figure}[h!]
\centerline{\includegraphics{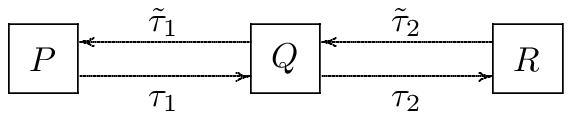}}
\end{figure}
}
\noindent
The set of tokens of the reduction of this medium to $\QQQ=\{P,R\}$ is empty. Thus this reduction is not a medium (Axiom [M2] is not satisfied).

The image $(\aa(\SSS),\bb(\TTT))$ of a token system $(\SSS,\TTT)$ under embedding $(\aa,\bb):(\SSS,\TTT)\rightarrow(\SSS',\TTT')$ is not, generally speaking, the reduction of $(\SSS',\TTT')$ to $\aa(\SSS)$. Indeed, let $\SSS'=\SSS$, and let $\TTT$ be a proper nonempty subset of $\TTT'$. Then the image of $(\SSS,\TTT)$ under the identity embedding is not the reduction of $(\SSS,\TTT')$ to $\SSS$ (which is $(\SSS,\TTT')$ itself).

On the other hand, this is true in the case of media as the following proposition demonstrates.

\begin{proposition}
Let $(\aa,\bb):(\SSS,\TTT)\rightarrow(\SSS',\TTT')$ be an embedding of a medium $(\SSS,\TTT)$ into a medium $(\SSS',\TTT')$. Then the reduction $(\aa(\SSS),\TTT'_{\aa(\SSS)})$ is isomorphic to $(\SSS,\TTT)$.
\end{proposition}

\begin{proof}
For $\tau\in\TTT$, we define $\bb'(\tau)=\bb(\tau)_{\aa(\SSS)}$, the reduction of $\bb(\tau)$ to $\aa(\SSS)$. Let $S\tau=T$ for $S\not=T$ in $\SSS$. Then $\aa(S)\bb(\tau)=\aa(T)$ for $\aa(S)\not=\aa(T)$ in $\aa(\SSS)$. Hence, $\bb'$ maps $\TTT$ to $\TTT'_{\aa(\SSS)}$.

Let us show that $(\aa,\bb')$ is an isomorphism from $(\SSS,\TTT)$ onto $(\aa(\SSS),\TTT'_{\aa(\SSS)})$.
\begin{itemize}
	\item[(i)] $\bb'$ is onto. Suppose $\tau'_{\aa(\SSS)}\not=\tau_0$ for some $\tau'\in\TTT'$. Then there are $P\not= Q$ in $\SSS$ such that $\aa(P)\tau'_{\aa(\SSS)}=\aa(P)\tau'=\aa(Q)$. Let $Q=P\bm$ where $\bm$ is a straight message. We have
$$
\aa(Q)=\aa(P\bm)=\aa(P)\bb(\bm)=\aa(P)\tau',
$$
implying, by Theorem~\ref{Theorem1.16}, $\bb(\bm)=\tau'$, since $\bb(\bm)$ is a straight message. Hence, $\bm=\tau$ for some $\tau\in\TTT$. Thus $\bb(\tau)=\tau'$, which implies
$$
\bb'(\tau)=\bb(\tau)_{\aa(\SSS)}=\tau'_{\aa(\SSS)}.
$$
	\item[(ii)] $\bb'$ is one--to--one. Suppose $\bb'(\tau_1)=\bb'(\tau_2)$. Since $\bb'(\tau_1)$ and $\bb'(\tau_2)$ are tokens on $\aa(\SSS)$ and $(\SSS',\TTT')$ is a medium, we have $\bb(\tau_1)=\bb(\tau_2)$. Hence, $\tau_1=\tau_2$.
	\item[(iii)]Finally,
$$
S\tau=T \quad \eq \quad \aa\left(S\right)\bb'\left(\tau\right)=\aa\left(T\right),
$$
since
$$
S\tau=T \quad \eq \quad \aa\left(S\right)\bb\left(\tau\right)=\aa\left(T\right).
$$
\end{itemize}
\end{proof}

We shall see later (Theorem~\ref{FiniteRepresentationTheorem}) that any medium is isomorphic to a submedium of a complete medium.

\section{Families of sets representable as media}

In this section, the objects of our study are token systems which are defined by means of families of subsets of a given set $X$.

Let $\POW(X)=2^X$ be the family of all subsets of $X$ and let $\GGG=\cup\{\gg_x,\tilde{\gg}_x\}$ be the family of functions from $\POW(X)$ into $\POW(X)$ defined by
\begin{align*}
&\gg_x:S\mapsto S\gg_x = S\cup\{x\}, \\
&\tilde{\gg}_x:S\mapsto S\tilde{\gg}_x = S\setminus\{x\}
\end{align*}
for all $x\in X$.
It is clear that $(\POW(X),\GGG)$ is a token system and that, for a given $x\in X$, the token $\tilde{\gg}_x$ is the unique reverse of the token $\gg_x$.

Let $\FFF$ be a nonempty family of subsets of $X$. In what follows, $(\FFF,\GGG_\FFF)$ stands for the reduction of the token system $(\POW(X),\GGG)$ to $\FFF\SB\POW(X)$.

In order to characterize token systems $(\FFF,\GGG_\FFF)$ which are media, we introduce some geometric concepts in $\POW(X)$ (cf. \cite{kB73,vKsO75,sO83}).

\begin{definition} \label{D:line segment}
Given $P,Q\in \POW(X)$, the \emph{interval} $[P,Q]$ is defined by 
\begin{equation*}
[P,Q]=\{R\in 2^X:P\cap Q\SB R\SB P\cup Q\}.
\end{equation*}
If $R\in[P,Q]$, we say that $R$ \emph{lies between $P$ and $Q$}. 

A sequence $P=P_0,P_1,\ldots,P_n=Q$ of distinct elements of $\POW(X)$ is a \emph{line segment} between $P$ and $Q$ if
\begin{enumerate}
	\item[\emph{L1.}] $P_i\in[P_k,P_m]$ for $k\leq i\leq m,$ and
	\item[\emph{L2.}] $R\in[P_i,P_{i+1}]$ implies $R=P_i$ or $R=P_{i+1}$ for all $0\leq i\leq n-1$.
\end{enumerate}

The \emph{distance}
between $P$ and $Q$ is defined by
\begin{equation*}
d(P,Q) = \begin{cases}
	|P\Delta Q|, &\text{if $P\Delta Q$ is a finite set,} \\
	\infty, &\text{otherwise,}
\end{cases}
\end{equation*}
where $\Delta$ stands for the symmetric difference operation.

A binary relation $\sim$ on $\POW(X)$ is defined by
\begin{equation*}
P\sim Q\quad\eq\quad d(P,Q)<\infty.
\end{equation*}
The relation $\sim$ is an equivalence relation on $\POW(X)$. We denote $[S]$ the equivalence class of $\sim$ containing $S\in\POW(X)$. We also denote $\POWF(X)=[\es]$, the family of all finite subsets of the set $X$.
\end{definition}

\begin{theorem} \label{DistanceTheorem}
Given $S\in \POW(X)$,
\begin{itemize}
	\item[\emph{(i)}] The distance function $d$ defines a metric on $[S]$.
	\item[\emph{(ii)}] $R$ lies between $P$ and $Q$ in $[S]$, that is, $R\in[P,Q],$ if and only if 
\begin{equation*}
d(P,R)+d(R,Q)=d(P,Q).
\end{equation*}
	\item[\emph{(iii)}] A sequence $P=P_0,P_1,\ldots,P_n=Q$ is a line segment between $P$ and $Q$ in $[S]$ if and only if $d(P,Q)=n$ and $d(P_i,P_{i+1})=1,\;0\leq i\leq n-1$.
\end{itemize}
\end{theorem}

\begin{proof}
(i) It is clear that $d(P,Q)\geq 0$ and $d(P,Q)=0$ if and only if $P=Q$, and that $d(P,Q)=d(Q,P)$ for all $P,Q\in\langle S\rangle$. 

It remains to verify the triangle inequality. Let $S_1,S_2,S_3$ be three sets in $[S]$. Since these sets belong to the same equivalence class of the relation $\sim$, the following six sets
$$
V_i = (S_j\cap S_k)\setminus S_i,\;U_i = S_i\setminus(S_j\cup S_k)\quad\text{for $\{i,j,k\}=\{1,2,3\}$,}
$$
are finite. It is not difficult to verify that
$$
S_i\Delta S_j=U_i\cup V_j\cup U_j\cup V_i,
$$
with disjoint sets in the right hand side of the equality. We have
\begin{align} \label{so1 triangle id}
&\phantom{=(}|S_i\Delta S_j|+|S_j\Delta S_k|\notag\\
&=(|U_i|+|V_j|+|U_j|+|V_i|)+(|U_j|+|V_k|+|U_k|+|V_j|)\\
&=|S_i\Delta S_k|+2(|U_j|+|V_j|), \notag
\end{align}
which implies the triangle inequality.

(ii) By~(\ref{so1 triangle id}),
$$
|S_i\Delta S_j|+|S_j\Delta S_k|=|S_i\Delta S_k|
$$
if and only if $U_j=\es$ and $V_j=\es$, or, equivalently, if and only if
$$
S_i\cap S_k\SB S_j\SB S_i\cup S_j.
$$

(iii) (Necessity.) By Condition L2 of Definition~\ref{D:line segment},
$$
d(P_i,P_{i+1})=|P_i\Delta P_{i+1|}|=1.
$$
Indeed, if there were $x\in P_i\setminus P_{i+1}$ and $y\in P_{i+1}\setminus P_i$, then the set $R=(P_i\setminus\{x\})\cup\{y\}$ would lie strictly between $P_i$ and $P_{i+1}$, a contradiction. By Condition L1 of Definition~\ref{D:line segment} and part (ii) of the theorem, we have
$$
d(P,Q)=1+d(P_1,Q)=\cdots=\underbrace{1+1+\cdots+1}_{n}=n.
$$

(Sufficiency.) Let $P=P_0,P_1,\ldots,P_n=Q$ be a sequence of sets such that $d(P_i,P_{i+1})=1$ and $d(P,Q)=n$. Condition L2 of Definition~\ref{D:line segment} is clearly satisfied. 

By the triangle inequality, we have
$$
d(P,P_i)\leq i,\;\;d(P_i,P_j)\leq j-i,\;\;d(P_j,Q)\leq n-j
$$
for $i<j$. Let us add these inequalities and use the triangle inequality again. We obtain
$$
n=d(P,Q)\leq d(P,P_i)+d(P_i,P_j)+d(P_j,Q)\leq n.
$$
It follows that $d(P_i,P_j)=j-i$ for all $i<j$. In particular, 
$$
d(P_k,P_i)+d(P_i,P_m)=d(P_k,P_m)\quad\text{for $k\leq i\leq m$.}
$$
By part (ii) of the theorem, $P_i\in[P_k,P_m]$ for $k\leq i\leq m$, which proves Condition L1 of Definition~\ref{D:line segment}.
\end{proof}

The concept of a line segment seems to be similar to the concept of a straight message. The following lemma validates this intuition.

\begin{lemma} \label{StraightMessage}
Let $(\FFF,\GGG_\FFF)$ be a token system and let $P$ and $Q$ be two distinct sets in $\FFF$. A message $\bm=\tau_1\tau_2\cdots\tau_n$ producing $Q$ from $P$ is
straight if and only if 
$$
P_0=P,\;P_1=P_0\tau_1,\;\ldots,\;P_n=P_{n-1}\tau_n=Q
$$
is a line segment between $P$ and $Q$.
\end{lemma}

\begin{proof}
(Necessity.) We use induction on $n=\ell(\bm)$. Let $\bm=\tau_1$. Since $\tau_1$ is either $\gg_x$ or $\tilde{\gg}_x$ for some $x\in X$ and effective, either $Q=P\cup\{x\}$ or $Q=P\setminus\{x\}$ and $Q\not= P$. Therefore, $d(P,Q)=1$, that is $\{P,Q\}$ is a line segment. 

Now, let us assume that the statement holds for all straight messages $\bn$ with $\ell(\bn)=n-1$ and let $\bm=\tau_1\tau_2\cdots\tau_n$ be a straight message producing $Q$ from $P$. Clearly, $\bm_1=\tau_2\cdots\tau_n$ is a straight message producing $Q$ from $P_1=P\tau_1$ and $\ell(\bm_1)=n-1$. 

Suppose that $\tau_1=\gg_x$ for some $x\in X$. Since $\bm$ is stepwise effective, $x\notin P$. Since $\bm$ is consistent, $x\in Q$. Therefore $P_1=P\tau_1=P\cup\{x\}\in [P,Q]$. Suppose that $\tau_1=\tilde{\gg}_x$ for some $x\in X$. Since $\bm$ is stepwise effective, $x\in P$. Since it is consistent, $x\notin Q$. Again, $P_1=P\setminus\{x\}\in[P,Q]$. In either case, $d(P,P_1)=1$. By Theorem~\ref{DistanceTheorem}(ii) and the induction hypothesis, $d(P,Q)=d(P,P_1)+d(P_1,Q)=n$. By the induction hypothesis and Theorem~\ref{DistanceTheorem}(iii), the sequence
\begin{equation*}
P_0=P,P_1=P_0\tau_1,\ldots,P_n=P_{n-1}\tau_n=Q
\end{equation*}
is a line segment between $P$ and $Q$.

(Sufficiency.) Let $P_0=P,P_1=P_0\tau_1,\ldots,P_n=P_{n-1}\tau_n=Q$ be a line segment between $P$ and $Q$ for some message $\bm=\tau_1\cdots\tau_n$. Clearly, $\bm$ is stepwise effective. To prove consistency, we use induction on $n$. The
statement is trivial for $n=1$. Suppose it holds for all messages of length less than $n$ and let $\bm$ be a message of length $n$. Suppose $\bm$ is inconsistent. By the induction hypothesis, this can occur only if either $\tau_1=\gg_x$, $\tau_n=\tilde{\gg}_x$ or $\tau_1=\tilde{\gg}_x$, $\tau_n=\gg_x$ for some $x\in X$. In the former case, $x\in P_1,\;x\notin P,\;x\notin Q$. In the latter case, $x\notin P_1,\;x\in P,\;x\in Q$. In both cases, $P_1\notin[P,Q]$, a
contradiction.
\end{proof}

The following theorem is an immediate consequence of Lemma~\ref{StraightMessage}.

\begin{theorem} \label{Fin[S]}
If $(\FFF,\GGG_\FFF)$ is a medium, then $\FFF\SB [S]$ for some $S\in\POW(X)$.
\end{theorem}

Clearly, the converse of this theorem is not true. To characterize those token systems $(\FFF,\GGG_\FFF)$ which are media, we use the concept of a well graded family of sets \cite{jDjF97}. (See also \cite{vKsO75,sO80} and \cite{sO83} where the same concept was introduced as a ``completeness condition''.)

\begin{definition}
A family $\FFF\SB\POW(X)$ is \emph{well graded} if for any two distinct sets $P$ and $Q$ in $\FFF$, there is a sequence of sets $P=R_0,R_1,\ldots,R_n=Q$ such that $d(R_{i-1},R_i)=1$ for $i=1,\ldots,n$ and $d(P,Q)=n$.
\end{definition}

In other words, $\FFF$ is a well graded family if for any two distinct elements $P,Q\in\FFF$ there is a line segment between $P$ and $Q$ in $\FFF$. 

\begin{theorem} \label{T:wg gamma}
Let $\FFF$ be a well graded family of subsets of some set $X$. Then $x\in X$ defines tokens $\gg_x,\tilde{\gg}_x\in\GGG_\FFF$ if and only if $x\in\cup\FFF\setminus\cap\FFF$.
\end{theorem}

\begin{proof}
It is clear that elements of $X$ that are not in $\cup\FFF\setminus\cap\FFF$ do not define tokens in $\GGG_\FFF$.

Suppose that $x\in\cup\FFF\setminus\cap\FFF$. Then there are sets $P$ and $Q$ in $\FFF$ such that $x\in Q\setminus P$. Let $R_0,R_1,\ldots,R_n$ be a line segment in $\FFF$ between $P$ and $Q$. Then there is $i$ such that $R_{i+1}=R_i\cup\{x\}$ and $x\notin R_i$. Therefore, $R_{i+1}=R_i\gg_x$ and $R_i=R_{i+1}\tilde{\gg}_x$, that is, $\gg_x,\tilde{\gg}_x\in\GGG_\FFF$.
\end{proof}

In the rest of the paper we assume that a well graded family $\FFF$ of subsets of $X$ defining a token system $(\FFF,\GGG_\FFF)$ satisfies the following conditions:
\begin{equation} \label{E:wg gamma}
\cap\FFF=\es\quad\text{and}\quad\cup\FFF=X.
\end{equation}

We have the following theorem (cf.~Theorem~4.2 in~\cite{sOaD00}).

\begin{theorem} \label{WellGradedMedia}
A token system $(\FFF,\GGG_\FFF)$ is a medium if and only if $\FFF$ is a well graded family of subsets of $X$.
\end{theorem}

\begin{proof}
(Necessity.) Suppose $(\FFF,\GGG_\FFF)$ is a medium. By axiom [M2*], for given $P,Q\in\FFF$, there is a straight message producing $Q$ from $P$. By Lemma~\ref{StraightMessage}, there is a line segment in $\FFF$ between $P$ and $Q$. Hence $\FFF$ is well graded.

(Sufficiency.) Let $\FFF$ be a well graded family of subsets of $X$. We need to show that the four axioms defining a medium are satisfied for $(\FFF,\GGG_\FFF)$.

[M1]. Clearly, $\gg_x$ and $\tilde{\gg}_x$ are unique mutual reverses of each other.

[M2*]. Follows immediately from Lemma~\ref{StraightMessage}.

[M3]. (Necessity.) Let $\bm$ be a message which is stepwise effective for $P\in\FFF$ and ineffective for this state, that is, $P\bm=P$.

Let $\tau$ be a token in $\bm$ such that $\tilde\tau\notin\bm$. If $\tau=\gg_x$ for some $x\in X$, then $x\notin P$ and $x\in P\bm$, since $\bm$ is stepwise effective for $P$ and $\tilde\gg_x=\tilde\tau\notin\bm$. We have a contradiction, since $P\bm=P$. In a similar way, we obtain a contradiction assuming that $\tau=\tilde\gg_x$. Thus, for each token $\tau$ in $\bm$, there is an appearance of the reverse token $\tilde\tau$ in $\bm$.

Because $\bm$ is stepwise effective, the appearances of tokens $\tau$ and $\tilde{\tau}$ in $\bm$ must alternate. Suppose that the sequence of appearances of $\tau$ and $\tilde{\tau}$ begins and ends with $\tau=\gg_x$ (the argument is similar if $\tau=\tilde{\gg}_x$). Since the message $\bm$ is stepwise effective for $P$ and ineffective for this state, we must have $x\notin P$ and $x\in P\bm=P$, a contradiction. It follows that $\bm$ is vacuous.

(Sufficiency.) Let $\bm$ be a vacuous message which is stepwise effective for some state $P$. Since $\bm$ is vacuous, the number of appearances of $\gg_x$ in $\bm$ is equal to the number of appearances of $\tilde\gg_x$ for any $x\in X$. Because $\bm$ is stepwise effective, the appearances of tokens $\gg_x$ and $\tilde{\gg}_x$ in $\bm$ must alternate. It follows that $x\in P$ if and only if $x\in P\bm$, that is $P\bm=P$. Thus the message $\bm$ is ineffective for $P$.

[M4]. Suppose two straight messages $\bm$ and $\bn$ produce $R$ from $P$ and $Q$, respectively, that is, $R=P\bm$ and $R=Q\bn$. Let us assume that $\bm$ and $\bn$ are not jointly consistent, that is, that $\boldsymbol{mn}$ is inconsistent. Then there are two mutually reverse tokens $\tau$ and $\tilde{\tau}$ in $\boldsymbol{mn}$. Since $\bm$ and $\bn$ are straight messages, we may assume, without loss of generality, that $\tau=\gg_x$ is in $\bm$ and $\tilde{\tau}=\tilde{\gg}_x$ is in $\bn$ for some $x\in X$. Since $\bm$ is straight, $x\in R$. Since $\bn$ is straight, $x\notin R$, a contradiction.
\end{proof}

Clearly, for a given $S\in\POW(X)$, $[S]$ is a well graded family of subsets of $X$. Hence, $([S],\GGG_{[S]})$ is a medium. It is easy to see that any such medium is a complete medium (see Definition~\ref{D:complete-media}). The converse is also true as the following theorem asserts. 

\begin{theorem} \label{CompleteMedium}
A medium in the form $(\FFF,\GGG_\FFF)$ is complete if and only if $\FFF=[S]$ for some $S\in\POW(X)$.
\end{theorem}

\begin{proof}
We need to prove necessity only. Suppose that $(\FFF,\GGG_\FFF)$ is a complete medium. By Theorem~\ref{Fin[S]}, $\FFF\SB[S]$ for some $S\in\POW(X)$. For a given $P\in[S]$, let $m=d(P,S)$. We prove that $P\in\FFF$ by induction on $m$.

Let $m=1$. Then either $P=S\cup\{x\}$ or $P=S\setminus\{x\}$ for some $x\in X$ and $P\not= S$. In the former case, $x\notin S$ implying that $\tilde{\gg}_x$ is not effective on $S$. By completeness, $S\gg_x= P$. Thus $P\in\FFF$. Similarly, if $P=S\setminus\{x\}$, then $x\in S$ and $S\tilde{\gg}_x= P$.

Suppose that $Q\in\FFF$ for all $Q\in[S]$ such that $d(S,Q)=m$ and let $P$ be an element in $[S]$ such that $d(S,P)=m+1$. Then there exists $R\in\FFF$ such that $d(S,R)=m$ and $d(R,P)=1$. Since $[R]=[S]$, it follows from the argument in the previous paragraph that $P\in\FFF$.
\end{proof}

The following theorem shows that all media in the form $([S],\GGG_{[S]})$ are isomorphic.

\begin{theorem} \label{IsomorphismTheorem}
For any $S',S\in\POW(X)$, the media $([S'],\GGG_{[S']})$ and $([S],\GGG_{[S]})$ are isomorphic.
\end{theorem}

\begin{proof}
It suffices to consider the case when $S'=\es$.

We define $\aa:\POWF(X)\rightarrow[S]$ and $\bb:\GGG_{\POWF(X)}\rightarrow\GGG_{[S]}$ by
\begin{align*}
&\aa(P)=P\Delta S, \\
&\bb(\tau)= \begin{cases}
	\tilde{\tau} &\text{if $\tau=\gg_x$ or $\tau=\tilde{\gg}_x$ for $x\in S$,} \\
	\tau &\text{if $\tau=\gg_x$ or $\tau=\tilde{\gg}_x$ for $x\notin S$.}
\end{cases}
\end{align*}
Clearly, $\aa$ and $\bb$ are bijections. To prove that $P\tau=Q$ implies $\aa(P)\bb(\tau)=\aa(Q)$, let us consider the following cases:

\begin{enumerate}
	\item $\tau=\gg_x,\;x\in S$.
		\begin{enumerate}
			\item $x\in P$. Then $Q=P$ and
$$\aa(P)\bb(\tau) = (P\Delta S)\setminus\{x\} = P\Delta S = Q\Delta S.$$
			\item $x\notin P$. Then $Q=P\cup\{x\}$ and
$$\aa(P)\bb(\tau) = (P\Delta S)\setminus\{x\} = P\Delta(S\setminus\{x\}) = Q\Delta S.$$
		\end{enumerate}
	\item $\tau=\gg_x,\;x\notin S$.
		\begin{enumerate}
			\item $x\in P$. Then $Q=P$ and
$$\aa(P)\bb(\tau) = (P\Delta S)\cup\{x\} = P\Delta S = Q\Delta S.$$
			\item $x\notin P$. Then $Q=P\cup\{x\}$ and
$$\aa(P)\bb(\tau) = (P\Delta S)\cup\{x\} =  Q\Delta S.$$
		\end{enumerate}
\end{enumerate}
A similar argument proves the theorem in the case when $\tau=\tilde{\gg}_x$. It is also easy to verify the converse implication: $\aa(P)\bb(\tau)=\aa(Q)\;\Rightarrow\; P\tau=Q$.
\end{proof}

We summarize the results of this section as follows:

\medskip
1.~A token system $(\FFF,\GGG_\FFF)$ is a medium if and only if $\FFF$ is a well graded family of subsets of $X$.

2.~A complete medium in the form $(\FFF,\GGG_\FFF)$ is $([S],\GGG_{[S]})$ for some $S\in\POW(X)$ and all such media are isomorphic.

3.~Any medium in the form $(\FFF,\GGG_\FFF)$ is a submedium of a complete medium and isomorphic to a submedium of the complete medium $(\POWF(X),\GGG_{\POWF(X)})$ of all finite subsets of $X$.

\section{The representation theorem}

In this section we show that any medium is isomorphic to a medium in the form $(\FFF,\GGG_\FFF)$ where $\FFF$ is a well graded family of finite subsets of some set $X$. In our construction we employ the concept of `orientation'~\cite{jFsO02}.

\begin{definition} \label{D:orientation}
An \emph{orientation} of a medium
$(\SSS,\TTT)$ is a partition of its set of tokens into two classes $\TTT^{+}$ and $\TTT^{-}$, respectively called \emph{positive} and \emph{negative}, such that for any $\tau\in\TTT$, we have 
$$
\tau\in\TTT^{+}\quad\eq\quad\tilde{\tau}\in\TTT^{-}
$$
A medium $\left(\SSS,\TTT\right)$ equipped with an orientation $\{\TTT^{+},\TTT^{-}\}$ is said to be \emph{oriented} by $\{\TTT^{+},\TTT^{-}\}$, and tokens from $\TTT^{+}$
\emph{(}resp. $\TTT^{-}$\emph{)} are called \emph{positive} \emph{(}resp. \emph{negative}\emph{)}. The \emph{positive content} \emph{(}resp.~\emph{negative content}\emph{)} of a state $S$ is the set $\widehat{S}^{+}=\widehat{S}\cap\TTT^{+}$ \emph{(}resp. $\widehat{S}^{-}=\widehat{S}\cap\TTT^{-}$\emph{)} of its positive \emph{(}resp. negative\emph{)} tokens.
\end{definition}

Let $(\SSS,\TTT)$ be a medium equipped with an orientation $\{\TTT^{+},\TTT^{-}\}$. In what follows, we show that $(\SSS,\TTT)$ is isomorphic to the medium $(\FFF,\GGG_\FFF)$ of the well graded family $\FFF=\{\widehat{S}^+\}_{S\in\SSS}$ of all positive contents.

\begin{lemma} \label{S+=T+}
\emph{(cf.~\cite{sOaD00})} For any two states $S,T\in\SSS$, 
\begin{equation*}
\widehat{S}^+=\widehat{T}^+\quad\eq\quad S=T.
\end{equation*}
\end{lemma}

\begin{proof}
By Theorem~\ref{Theorem1.17-2}, it suffices to prove that $\widehat{S}^+=\widehat{T}^+$ implies $\widehat{S}=\widehat{T}$. Let $\tau\in\widehat{S}^-$. Then, by Theorem~\ref{Theorem1.17-1}, $\tilde{\tau}\notin\widehat{S}^+$. Hence, $\tilde{\tau}\notin\widehat{T}^+$ which implies $\tau\in\widehat{T}^-$. Therefore $\widehat{S}^-\SB\widehat{T}^-$. By symmetry, $\widehat{S}^-=\widehat{T}^-$. Hence, $\widehat{S}=\widehat{T}$.
\end{proof}

We define
\begin{equation} \label{alpha}
\aa : S\mapsto \widehat{S}^+,
\end{equation}
for $S\in\SSS$. It follows from Lemma~\ref{S+=T+} that $\aa$ is a bijection.

Suppose that $\tau\in\cap\,\FFF=\cap_{S\in\SSS}\widehat{S}^+$. There are $S,T\in\SSS$ such that $T=S\tau$. Then, by Theorem~\ref{Theorem1.16}, $\widehat{T}\setminus\widehat{S}=\{\tau\}$, that is, $\tau\notin\widehat{S}\supseteq\widehat{S}^+$. Hence, $\cap\,\FFF=\es$.

Let $\tau\in\TTT^+$. There are $S,T\in\SSS$ such that $T=S\tau$. Then $\tau\in\widehat{T}^+$. Hence, $\cup\,\FFF=\TTT^+$.

We define a mapping $\bb:\TTT\rightarrow\GGG_\FFF$ by
\begin{equation} \label{beta}
\bb(\tau)=\begin{cases}
	\gg_\tau &\text{if $\tau\in\TTT^+$,} \\
	\tilde{\gg}_{\tilde{\tau}} &\text{if $\tau\in\TTT^-$,}
\end{cases}
\end{equation}
and show that the pair $(\aa,\bb)$, where $\aa$ and $\bb$ are mappings defined, respectively, by (\ref{alpha}) and (\ref{beta}), is an isomorphism between the token systems $(\SSS,\TTT)$ and $(\FFF,\GGG_\FFF)$.

\begin{theorem} \label{RepresentationTheorem}
Let $(\SSS,\TTT)$ be an oriented medium. For all $S,T\in\SSS$ and $\tau\in\TTT$,
\begin{equation*}
T=S\tau\quad\eq\quad\aa(T)=\aa(S)\bb(\tau),
\end{equation*}
that is, the token systems $(\SSS,\TTT)$ and $(\FFF,\GGG_\FFF)$ are isomorphic.
\end{theorem}

\begin{proof}
(i) Suppose that $\tau\in\TTT^+$. We need to prove that
\begin{equation*}
T=S\tau\quad\eq\quad\widehat{T}^+=\widehat{S}^+\gg_\tau
\end{equation*}
for all $S,T\in\SSS$ and $\tau\in\TTT$.

Let us consider the following cases:
\begin{enumerate}
	\item $\tau\in\widehat{S}^+$. Suppose $S\tau=T\not=S$. Then, by Theorem~\ref{Theorem1.16}, $\tau\in\widehat{T}\setminus\widehat{S}$, a contradiction. Hence, $T=S$. Clearly, 
$\widehat{S}^+=\widehat{S}^+\cup\{\tau\}=\widehat{S}^+\gg_\tau$.
	\item $\tau\notin\widehat{S}^+,\,\widehat{S}^+\cup\{\tau\}\notin\FFF$. Suppose that $S\tau=T\not=S$. By Theorem~\ref{Theorem1.16}, $\widehat{T}\setminus\widehat{S}=\{\tau\}$ and $\widehat{S}\setminus\widehat{T}=\{\tilde{\tau}\}$. Since $\tau\in\TTT^+$, we have $\widehat{S}^+\cup\{\tau\}=\widehat{T}^+\in\FFF$, a contradiction. Hence, in this case, $S=S\tau$ and $\widehat{S}^+=\widehat{S}^+\gg_\tau$.
	\item $\tau\notin\widehat{S}^+,\,\widehat{S}^+\cup\{\tau\}\in\FFF$. Then there exists $T\in\SSS$ such that $\widehat{T}^+=\widehat{S}^+\cup\{\tau\}$. Thus $\tau\in\widehat{T}\setminus\widehat{S}$. Suppose that there is $\tau'\not=\tau$ which is also in $\widehat{T}\setminus\widehat{S}$. Then $\tau'$ is a negative token. Since $\tau'\notin\widehat{S}$, we have $\tilde{\tau}'\in\widehat{S}^+\subset\widehat{T}^+\SB\widehat{T}$. Hence, $\tau'\notin\widehat{T}$, a contradiction. It follows that $\widehat{T}\setminus\widehat{S}=\{\tau\}$. By Theorem~\ref{Theorem1.16}, $T=S\tau$. By the argument in item 2, $\widehat{T}^+=\widehat{S}^+\cup\{\tau\}=\widehat{S}^+\gg_\tau$.
\end{enumerate}

(ii) Suppose that $\tau\in\TTT^-$. Then $\tilde{\tau}\in\TTT^+$. We need to prove that
$$
T=S\tau\quad \eq\quad\widehat{T}^+=\widehat{S}^+\tilde{\gg}_{\tilde{\tau}}
$$
for all $S,T\in\SSS$ and $\tau\in\TTT$.

Let us consider the following cases:
\begin{enumerate}
	\item $\tilde{\tau}\notin\widehat{S}^+$. Suppose that $S\tau=T\not=S$. Then $S=T\tilde{\tau}$ and, by Theorem~\ref{Theorem1.16}, $\tilde{\tau}\in\widehat{S}\setminus\widehat{T}$, a contradiction since $\tilde{\tau}$ is a positive token. Hence, $S=S\tau$. On the other hand, $\widehat{S}^+=\widehat{S}^+\setminus\{\tilde{\tau}\}=\widehat{S}^+\tilde{\gg}_{\tilde{\tau}}$.
	\item $\tilde{\tau}\in\widehat{S}^+,\,\widehat{S}^+\setminus\{\tilde{\tau}\}\notin\FFF$. Suppose again that $S\tau=T\not=S$. Then $S=T\tilde{\tau}$ and, by Theorem~\ref{Theorem1.16}, $\{\tilde{\tau}\}=\widehat{S}\setminus\widehat{T}$ and $\{\tau\}=\widehat{T}\setminus\widehat{S}$. Since $\tilde{\tau}$ is a positive token, we have $\widehat{S}^+\setminus\{\tilde{\tau}\}=\widehat{T}^+$, a contradiction. Hence, in this case, $S=S\tau$ and $\widehat{S}^+=\widehat{S}^+\tilde{\gg}_{\tilde{\tau}}$.
	\item $\tilde{\tau}\in\widehat{S}^+,\,\widehat{S}^+\setminus\{\tilde{\tau}\}\in\FFF$. There is $T\in\SSS$ such that $\widehat{T}^+=\widehat{S}^+\setminus\{\tilde{\tau}\}$. We have $\tau\notin\widehat{S}$, since $\tilde{\tau}\in\widehat{S}^+$, and $\tau\in\widehat{T}$, since $\tilde{\tau}\notin\widehat{T}^+$ and $\tilde{\tau}$ is a positive token. Hence, $\tau\in\widehat{T}\setminus\widehat{S}$. Suppose that there is $\tau'\not=\tau$ which is also in $\widehat{T}\setminus\widehat{S}$. Then $\tau'$ is a negative token. Since $\tau'\notin\widehat{S}$, we have $\tilde{\tau}'\in\widehat{S}^+$. Since $\tau'\not=\tau$, we have $\tilde{\tau}'\in\widehat{T}^+$. Hence, $\tau'\notin\widehat{T}$, a contradiction. It follows that $\widehat{T}\setminus\widehat{S}=\{\tau\}$. By Theorem~\ref{Theorem1.16}, $T=S\tau$. Clearly, $\widehat{T}^+=\widehat{S}^+\setminus\{\tau\}=\widehat{S}^+\tilde{\gg}_{\tilde{\tau}}$.
\end{enumerate}
\end{proof}

Since $(\SSS,\TTT)$ is a medium, the token system $(\FFF,\GGG_\FFF)$ is also a medium. By Theorem~\ref{WellGradedMedia}, we have the following result.

\begin{corollary}
$\FFF=\{\widehat{S}^+\}_{S\in\SSS}$ is a well graded family of subsets of $\TTT^+$.
\end{corollary}

Theorem~\ref{RepresentationTheorem} states that an oriented medium is isomorphic to the medium of its family of positive contents. By `forgetting' the orientation, one can say that any medium $(\SSS,\TTT)$ is isomorphic to some medium $(\FFF,\GGG_\FFF)$ of a well graded family $\FFF$ of sets. The following theorem is a stronger version of this result (cf.~Theorem~\ref{IsomorphismTheorem}).

\begin{theorem} \label{FiniteRepresentationTheorem}
Any medium $(\SSS,\TTT)$ is isomorphic to a medium of a well graded family of finite subsets of some set $X$.
\end{theorem}

\begin{proof}
Let $S_0$ be a fixed state in $\SSS$. By Theorem~\ref{Theorem1.17-1}, the state $S_0$ defines an orientation $\{\TTT^{+},\TTT^{-}\}$ with $\TTT^-=\widehat{S}_0$ and $\TTT^{+}=\TTT\setminus\TTT^-$. Let $\bm$ be a straight message producing a state $S$ from the state $S_0$. By Theorem~\ref{Theorem1.16}, we have $\CCC(\bm)=\widehat{S}\setminus\widehat{S}_0=\widehat{S}^+$. Thus, $\widehat{S}^+$ is a finite set. The statement of the theorem follows from Theorem~\ref{RepresentationTheorem}.
\end{proof}

\begin{remark}
{\rm An infinite oriented medium may have infinite positive contents of all its states. Consider, for instance, the medium $(\FFF,\GGG_\FFF)$, where 
$$
\FFF=\{S_n=]-\infty,n]:n\in\Zee\}.
$$
Then each $\widehat{S}_n^+=\{\gg_k\}_{k<n}$ is an infinite set. Nevertheless, by Theorem~\ref{FiniteRepresentationTheorem}, the medium $(\FFF,\GGG_\FFF)$ is isomorphic to a medium of a well graded family of finite sets.
}
\end{remark}

\section{Media and graphs}

In this section we study connections between media and graph theories.

\begin{definition}
Let $(\SSS,\TTT)$ be a medium. We say that a graph $G=(V,E)$ \emph{represents} $(\SSS,\TTT)$ if there is a bijection $\aa:\SSS\rightarrow V$ such that an unordered pair of vertices $PQ$ is an edge of the graph if and only if $P\not=Q$ (no loops) and there is $\tau\in\TTT$ such that $\aa^{-1}(P)\tau=\aa^{-1}(Q)$.

A graph $G$ representing a medium $(\SSS,\TTT)$ is the \emph{graph of the medium} $(\SSS,\TTT)$ if the set of states $\SSS$ is the set of vertices of $G$, the mapping $\aa$ is the identity, and the edges of $G$ are defined as above.
\end{definition}

Clearly, any graph which is isomorphic to a graph representing a token system $(\SSS,\TTT)$, also represents $(\SSS,\TTT)$, and isomorphic media are represented by isomorphic graphs.

The main goal of this section is to show that two media which are represented by isomorphic graphs are isomorphic (Theorem~\ref{Medium=Graph}).

First, we prove three lemmas about media and their graph representations.

\begin{lemma} \label{Lemma6.1}
Let $(\mathcal{S,T})$ be a medium, and suppose that: $\tau\in\TTT$ is a token; $S,T,P$ and $Q$ are states in $\SSS$ such that $S\tau=T$ and $P\tau=Q$. Let $\bm$ and $\bn$ be two straight messages producing $P$ from $S$ and $Q$ from $T$, respectively. Then $\bm$ and $\bn$ have equal contents and lengths, that is  $\CCC(\bm)=\CCC(\bn)$ and $\ell(\bm)=\ell(\bn)$, and $\bm\tau$ and $\tau\bn$ are straight messages for $S$.
\end{lemma}

\begin{proof}
The message $\bm\tau\tilde{\bn}\tilde{\tau}$ is stepwise effective for $S$ and ineffective for that state. By Axiom [M3], this message is vacuous. Hence, $\CCC(\bm)=\CCC(\bn)$ and $\ell(\bm)=\ell(\bn)$.
Each of two straight messages $\bm$ and $\tilde{\tau}$ produces $P$. By Axiom [M4], they are jointly consistent, that is, $\tau\notin\CCC(\bm)$. Hence, $\bm\tau$ is a straight message. Similarly, $\tau\bn$ is a straight message.
\end{proof}

\begin{lemma} \label{Lemma6.2}
Let $S,T,P,Q$ be four distinct states such that $S\tau_1=T$, $P\tau_2=Q$ and $\bm$ and $\bn$ be two straight messages producing $P$ from $S$ and $Q$ from $T$, respectively. If the messages $\tau_1\bn,\bm\tau_2$, and $\tilde{\tau}_1\bm$ are straight, then $\tau_1=\tau_2$.
\end{lemma}

\begin{proof}
Suppose that $\tau_1\not=\tau_2$. By Theorem~\ref{Theorem1.14}, $\CCC(\tau_1\bn)=\CCC(\bm\tau_2)$. Hence, $\tau_1\in\CCC(\bm)$, a contradiction, since we assumed that $\tilde{\tau}_1\bm$ is a straight message.
\end{proof}

\begin{lemma} \label{Lemma6.3}
Let $(\SSS,\TTT)$ be a medium, $G=(V,E)$ be a graph representing this medium, and $\aa$ be the bijection $\SSS\rightarrow V$ defining the graph $G$. If $\bm=\tau_1\cdots\tau_m$ is a straight message transforming a state $S$ into a state $T$, then the sequence of vertices $(\aa(S_i))_{0\leq i\leq m}$, where $S_i=S\tau_0\tau_1\cdots\tau_i$, forms a shortest path joining $\aa(S)$ and $\aa(T)$ in $G$. Conversely, if a sequence $(\aa(S_i))_{0\leq i\leq m}$ is a shortest path connecting $\aa(S_0)=\aa(S)$ and $\aa(S_m)=\aa(T)$, then $S\bm= T$ for some straight message $\bm$ of length $m$.
\end{lemma}

\begin{proof}
(Necessity.) Let $\aa(P_0)=\aa(S),\aa(P_1),\ldots,\aa(P_n)=\aa(T)$ be a path in $G$ joining $\aa(S)$ and $\aa(T)$. There is a stepwise effective message $\bn=\rho_1\cdots\rho_n$ such that $P_i=T\rho_1\cdots\rho_{n-i}$ for $0\leq i<n$. The message $\boldsymbol{mn}$ is stepwise effective for $S$ and ineffective for this state. By Axiom [M3], this message is vacuous. Since $\bm$ is a straight message for $S$, we have $\ell(\bm)\leq\ell(\bn)$. It follows that $(\aa(S_i))_{0\leq i\leq m}$ forms a shortest path joining $\aa(S)$ and $\aa(T)$ in $G$.

(Sufficiency.) Let $\aa(S_0)=\aa(S),\aa(S_1),\ldots,\aa(S_m)=\aa(T)$ be a shortest path connecting vertices $\aa(S)$ and $\aa(T)$ in $G$. Then $S_i\tau_{i+1}=S_{i+1}$ for some tokens $\tau_i$, $1\leq i\leq m$. The message $\bm=\tau_1\cdots\tau_m$ transforms the state $S$ into the state $T$. By the argument in the necessity part of the proof, $\bm$ is a straight message for $S$.
\end{proof}

\begin{definition}\label{so1 cube and partial cube}
{\rm (cf.~\cite{dD73,wI00})} Let $X$ be a set. The graph $\HHH(X)$ is defined as follows: the set of vertices is the set $\POWF(X)$ of all finite subsets of $X$; two vertices $P$ and $Q$ are adjacent if the symmetric difference $P\Delta Q$ is a singleton. We say that $\HHH(X)$ is a \emph{cube on} $X$. Isometric subgraphs of the cube $\HHH(X)$, as well as graphs that are isometrically embeddable in $\HHH(X)$, are called \emph{partial cubes}.
\end{definition}

The proof of the following proposition is straightforward and omitted.

\begin{proposition} \label{cube=wg-family}
An induced subgraph $G=(V,E)$ of the cube $\HHH(X)$ is a partial cube if and only if $V$ is a well graded family of finite subsets of $X$. Then a shortest path in $G$ is a line segment in $\HHH(X)$ and the graph distance function $d$ on both $\HHH(X)$ and $G$ is given by
$$
d(P,Q)=|P\Delta Q|.
$$
\end{proposition}

The following theorem characterizes media in terms of their graphs.

\begin{theorem} \label{MediumGraph}
A graph $G$ represents a medium $(\SSS,\TTT)$ if and only if $G$ is a partial cube.
\end{theorem}

We give two proofs of this important theorem. The first proof uses the representation theorem.

\begin{proof}
(Necessity.) By Theorem~\ref{FiniteRepresentationTheorem}, we may assume that the given medium is $(\FFF,\GGG_\FFF)$ where $\FFF$ is a well graded family of finite subsets of some set $X$. By Proposition~\ref{cube=wg-family}, the graph of this medium is a partial cube.

(Sufficiency.) For a partial cube $G=(V,E)$ there is an isometric embedding $\aa$ of $G$ into a cube $\mathcal{H}(X)$ for some set $X$. Vertices of $\aa(G)$ form a well graded family of subsets of $X$. Then the medium $(\aa(V),\GGG_{\aa(V)})$ has $G$ as its graph.
\end{proof}

The second proof utilizes the Djokovi\'{c}--Winkler relation $\Theta$ on the set of edges of a graph. The definition and properties of this relation are found in the book~\cite{wI00}.

\begin{proof}
(Necessity.) We may assume that $G$ is the graph of $(\SSS,\TTT)$. Let 
$$
S,S_1,\ldots,S_n=S
$$
be a cycle of length $n$ in $G$. There is a stepwise effective message $\bm$ such that $S\bm=S$. By Axiom [M3], $\bm$ is vacuous. Therefore, $\ell(\bm)=n$ is an even number. Hence, $G$ is a bipartite graph.

For each edge $ST$ of $G$ there is a unique unordered pair of tokens $\{\tau,\tilde{\tau}\}$ such that $S\tau=T$ and $T\tilde{\tau}=S$. We denote $\sim$ the equivalence relation on the set of edges of $G$ defined by this correspondence. Let $ST\sim PQ$ and the notation is chosen such that $S\tau=T$ and $P\tau=Q$. Then, by lemmas~\ref{Lemma6.1} and~\ref{Lemma6.3},
\begin{equation} \label{Eq6.1}
d=d(S,P)=d(T,Q)=d(S,Q)-1=d(T,P)-1
\end{equation}

By Lemma~2.3 in~\cite{wI00}, $ST\Theta PQ$. On the other hand, if $ST\Theta PQ$ holds, then, by the same lemma, equation~(\ref{Eq6.1}) holds. By lemmas~\ref{Lemma6.2} and~\ref{Lemma6.3}, there is a token $\tau$ such that $S\tau=T$ and $P\tau=Q$, that is, $ST\sim PQ$. Thus $\sim\,=\Theta$, that is, $\Theta$ is an equivalence relation. By Theorem~2.10 in~\cite{wI00}, $G$ is a partial cube.

(Sufficiency.) We have already shown in the first proof that a partial cube is the graph of a medium.
\end{proof}

\begin{theorem} \label{Medium=Graph}
Two media are isomorphic if and only if the graphs representing these media are isomorphic.
\end{theorem}

\begin{proof}
(Necessity.) Clearly, graphs representing isomorphic media are isomorphic.

(Sufficiency.) Let $(\SSS,\TTT)$ and $(\SSS',\TTT')$ be two media and $G=(V,E)$ and $G'=(V',E')$ be two isomorphic graphs representing these media. Since $G$ and $G'$ are isomorphic, $G$ represents $(\SSS',\TTT')$. Thus we need to show that two media represented by the same graph are isomorphic.

Since $G$ is a partial cube, it also represents a medium $(\FFF,\GGG_\FFF)$ of a well graded family of subsets of some set $X$. We denote $\mu:\SSS\rightarrow V$ and $\nu:V\rightarrow\FFF$ the two bijections that define the graph representation $G$ of $(\SSS,\TTT)$ and $(\FFF,\GGG_\FFF)$, respectively. Then $\aa=\nu\circ\mu$ is a bijection $\SSS\rightarrow\FFF$ such that 
$$
S\tau = T\quad\eq\quad |\aa(S)\Delta\aa(T)|=1,
$$
for all $S\not= T$ in $\SSS$ and $\tau\in\TTT$.

Clearly, it suffices to prove that the media $(\SSS,\TTT)$ and $(\FFF,\GGG_\FFF)$ are isomorphic.

Let $\tau$ be a token in $\TTT$ and $S$ and $T$ be two distinct states in $\SSS$ such that $S\tau= T$. Then either $\aa(T)=\aa(S)\cup\{x\}$ for some $x\notin \aa(S)$ or $\aa(T)=\aa(S)\setminus\{x\}$ for some $x\in \aa(S)$. We define $\bb:\TTT\rightarrow\GGG_\FFF$ by
$$
\bb(\tau) = \begin{cases}
	\gg_x, &\text{if $\aa(T)=\aa(S)\cup\{x\}$ for some $x\notin\aa(S)$,} \\
	\tilde{\gg}_x, &\text{if $\aa(T)=\aa(S)\setminus\{x\}$ for some $x\in\aa(S)$.}
\end{cases}
$$
Let us show that $\bb$ does not depend on the choice of $S$ and $T$. We consider only the case when $\bb(\tau)=\tau_x$. The other case is treated similarly.

Let $P,Q$ be another pair of distinct states in $\SSS$ such that $P\tau= Q$, and let $P=S\bm$ and $Q=T\bn$ for some straight messages $\bm$ and $\bn$. By Lemma~\ref{Lemma6.1}, $\ell(\bm)=\ell(\bn)$. Then, by Lemma~\ref{Lemma6.3}, $d(\aa(S),\aa(P))=d(\aa(T),\aa(Q))$, and, by Lemma~\ref{Lemma6.1},
\begin{align*}
d(\aa(S),\aa(Q))&=d(\aa(S),\aa(T))+d(\aa(T),\aa(Q)),\\
d(\aa(T),\aa(P))&=d(\aa(T),\aa(S))+d(\aa(S),\aa(P)).
\end{align*}
By Theorem~\ref{DistanceTheorem},
\begin{gather*}
\aa(S)\cap\aa(Q)\;\SB\:\aa(T)=\aa(S)\cup\{x\}\;\SB\;\aa(S)\cup\aa(Q),\\
\aa(T)\cap\aa(P)=[\aa(S)\cup\{x\}]\cap\aa(P)\;\SB\;\aa(S)\;\SB\;\aa(S)\cup\aa(P).
\end{gather*}
Since $x\notin\aa(S)$, it follows that $x\in\aa(Q)$ and $x\notin\aa(P)$. Then 
\begin{equation*}
\aa(Q)=\aa(P)\cup\{x\},
\end{equation*}
since $d(\aa(P),\aa(Q))=1$. Hence, the mapping $\bb:\SSS\rightarrow\GGG_\FFF$ is well defined. 

Clearly, $\bb$ is a bijection satisfying the condition
$$
S\tau= T\quad\eq\quad\aa(S)\bb(\tau)=\aa(T).
$$
Therefore $(\aa,\bb)$ is an isomorphism from $(\SSS,\TTT)$ onto $(\FFF,\GGG_\FFF)$.
\end{proof}

We conclude this section with an example illustrating Theorem~\ref{Medium=Graph}.

\begin{example}
{\rm If $(\mathcal{S,T})$ and $(\mathcal{S',T'})$ are two finite isomorphic media, then $|\SSS|=|\mathcal{S'}|$ and $|\TTT|=|\mathcal{T'}|$. The converse, generally speaking, is not true. Consider, for instance, two media, $(\FFF,\GGG_\FFF)$ and $(\FFF',\GGG_{\FFF'})$, of well graded subsets of $X=\{a,b,c\}$ with
\begin{equation*}
\FFF=\{a,b,c,ab,ac,bc\}\quad\text{and}\quad\FFF'=\{a,c,ab,ac,bc,abc\}.
\end{equation*}
Their graphs, $G$ and $G'$, are shown in Figure~\ref{G and G'}.

{\begin{figure}[h!]
\centerline{\includegraphics{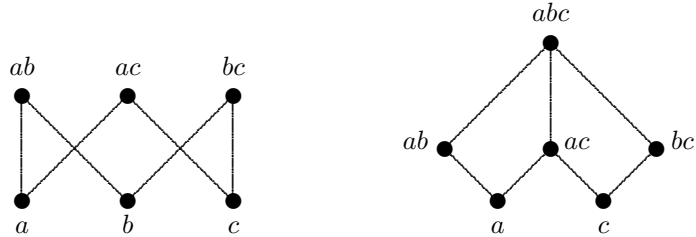}}
\caption{Graphs $G$ and $G'$.} \label{G and G'}
\end{figure}
}

Clearly, these graphs are not isomorphic. Thus the media $(\FFF,\GGG_\FFF)$ and $(\FFF',\GGG_{\FFF'})$ are not isomorphic.
}
\end{example}

\section{Uniqueness of media representations}

Theorem~\ref{FiniteRepresentationTheorem} asserts that any medium $(\SSS,\TTT)$ is isomorphic to the medium $(\FFF,\GGG_\FFF)$ of a well graded family $\FFF$ of finite subsets of some set $X$. In this section we show that this representation is unique in some precise sense.

Let $(\FFF_1,\GGG_{\FFF_1})$ and $(\FFF_2,\GGG_{\FFF_2})$ be two isomorphic representations of $(\SSS,\TTT)$ with well graded families $\FFF_1$ and $\FFF_2$ of subsets of $X_1$ and $X_2$, respectively. By Theorem~\ref{T:wg gamma},
$$
|\TTT|=|\GGG_{\FFF_i}|=2|\cup\FFF_i\setminus\cap\FFF_i|\quad\text{for $i=1,2$}.
$$
Thus, without loss of generality, we may assume that $\FFF_1$ and $\FFF_2$ are well graded families of finite subsets of the same set $X$ and that they satisfy conditions~(\ref{E:wg gamma}). The graphs of the media $(\FFF_1,\GGG_{\FFF_1})$ and $(\FFF_2,\GGG_{\FFF_2})$ are isomorphic partial subcubes of the cube $\HHH(X)$. On the other hand, by theorems~\ref{Medium=Graph} and~\ref{MediumGraph}, isometric partial cubes represent isomorphic media.

We formulate the uniqueness problem geometrically as follows:
\begin{quote}
Show that any isometry between two partial subcubes of $\HHH(X)$ can be extended to an isometry of the cube $\HHH(X)$.
\end{quote}

In other words, we want to show that partial subcubes of $\HHH(X)$ are unique up to isometries of $\HHH(X)$ onto itself.

\begin{remark}
{\rm Note, that $\HHH(X)$ is not a fully homogeneous space~(as defined, for instance, in \cite{dB01}), that is, an isometry between two arbitrary subsets of $\HHH(X)$, generally speaking, cannot be extended to an isometry of the cube $\HHH(X)$. On the other hand, $\HHH(X)$ is a homogeneous metric space.
}
\end{remark}

A general remark is in order. Let $Y$ be a homogeneous metric space, $A$ and $B$ be two metric subspaces of $Y$, and $\aa$ be an isometry from $A$ onto $B$. Let $c$ be a fixed point in $Y$. For a given $a\in A$, let $b=\aa(a)\in B$. Since $Y$ is homogeneous, there are isometries $\bb$ and $\gg$ of $Y$ such that $\bb(a)=c$ and $\gg(b)=c$. Then $\ll=\gg\aa\bb^{-1}$ is an isometry from $\bb(A)$ onto $\gg(B)$ such that $\ll(c)=c$. Clearly, $\aa$ is extendable to an isometry of $Y$ if and only if $\ll$ is extendable. Therefore, in the case of the space $\HHH(X)$, we may consider only well graded families of subsets containing the empty set $\es$ and isometries between these families fixing this point.

In what follows, we assume that $\es\in\FFF$ and $\cup\FFF=X$.

\begin{definition}
We define
$$
r_\FFF(x)=\min\{|A|: x\in A, A\in\FFF\}
$$
and, for $k\geq 1$,
$$
X_k^\FFF=\{x\in X: r_\FFF(x)=k\}.
$$
\end{definition}

We have $X_i^\FFF\cap X_j^\FFF=\es$ for $i\not=j$, and $\cup_k X_k^\FFF=X$. Note that some of the sets $X_k^\FFF$ could be empty for $k>1$, although $X_1^\FFF$ is not empty, since, by the wellgradedness property, $\FFF$ contains at least one singleton (we assumed that $\es\in\FFF$).

\begin{example}
{\rm Let $X=\{a,b,c\}$ and
$$
\FFF=\{\es,\{a\},\{b\},\{a,b\},\{a,b,c\}\}.
$$
Clearly, $\FFF$ is well graded. We have $r_\FFF(a)=r_\FFF(b)=1,\;r_\FFF(c)=3$, and
$$
X_1^\FFF=\{a,b\},\;X_2^\FFF=\es,\;X_3^\FFF=\{c\}.
$$
}
\end{example}

\begin{lemma}
For $A\in\FFF$ and $x\in A$, we have
\begin{equation} \label{OneLess}
r_\FFF(x)=|A|\quad\Rightarrow\quad A\setminus\{x\}\in\FFF.
\end{equation}
\end{lemma}

\begin{proof}
Let $k=|A|$. Since $\FFF$ is well graded, there is a nested sequence $\{A_i\}_{0\leq i\leq k}$ of distinct sets in $\FFF$ with $A_0=\es$ and $A_k=A$. Since $r_\FFF(x)=k$, we have $x\notin A_i$ for $i<k$. Hence, $A\setminus \{x\}=A_{k-1}\in\FFF$.
\end{proof}

Let us recall (Theorem~\ref{DistanceTheorem}(ii)) that 
\begin{equation} \label{so1 betweenness}
B\cap C\SB A\SB B\cup C\;\eq\; d(B,A)+d(A,C)=d(B,C)
\end{equation}
for all $A,B,C\in\POW(X)$.

It follows that isometries between two well graded families of sets $\FFF_1$ and $\FFF_2$ preserve the betweenness relation, that is,
\begin{equation} \label{so1 inclusion}
B\cap C\SB A\SB B\cup C\;\eq\;\aa(B)\cap \aa(C)\SB \aa(A)\SB \aa(B)\cup \aa(C)
\end{equation}
for $A,B,C\in\FFF_1$ and an isometry $\aa:\FFF_1\rightarrow\FFF_2$.

In the sequel, $\FFF_1$ and $\FFF_2$ are two well graded families of finite subsets of $X$ and $\aa:\FFF_1\rightarrow\FFF_2$ is an isometry such that $\aa(\es)=\es$. 

\begin{definition}
We define a binary relation $\pi$ between on $X$ by means of the following construction. By~(\ref{OneLess}), for a given $x\in X$ there is $A\in\FFF_1$ such that $x\in A$, $r_{\FFF_1}(x)=|A|$, and $A\setminus\{x\}\in\FFF_1$.
Since $\es\SB A\setminus\{x\}\subset A$, we have, by~(\ref{so1 inclusion}), $\aa(A\setminus\{x\})\subset\aa(A)$. Since $d\,(A\setminus\{x\},A)=1$, there is $y\in X,\;y\notin\aa(A)$ such that $\aa(A)=\aa(A\setminus\{x\})\cup\{y\}$. In this case we say that $xy\in\pi$.
\end{definition}

\begin{lemma}
If $x\in X_k^{\FFF_1}$ and $xy\in\pi$, then $y\in X_k^{\FFF_2}$.
\end{lemma}

\begin{proof}
Let $A\in\FFF_1$ be a set of cardinality $k$ defining $r_{\FFF_1}(x)=k$.
Since 
$$
|A|=d\,(\es,A)=d\,(\es,\aa(A))=|\aa(A)|\quad\text{and}\quad y\in\aa(A),
$$
we have $r_{\FFF_2}(y)\leq k$. Suppose that $m=r_{\FFF_2}(y)<k$. Then, by~(\ref{OneLess}), there is $B\in\FFF_2$ such that $y\in B$, $|B|=m$, and $B\setminus\{y\}\in\FFF_2$. Clearly,
$$
\aa(A\setminus\{x\})\cap B\SB\aa(A)\SB\aa(A\setminus\{x\})\cup B.
$$
By~(\ref{so1 inclusion}), we have
$$
(A\setminus\{x\})\cap\aa^{-1}(B)\SB A\SB (A\setminus\{x\})\cup\aa^{-1}(B).
$$
Thus, $x\in\aa^{-1}(B)$, a contradiction, since $r_{\FFF_1}(x)=k$ and $|\aa^{-1}(B)|=m<k$. It follows that $r_{\FFF_2}(y)=k$, that is, $y\in X_k^{\FFF_2}$.
\end{proof}

We proved that, for every $k\geq 1$, the restriction of $\pi$ to $X_k^{\FFF_1}$ is a relation $\pi_k$ between $X_k^{\FFF_1}$ and $X_k^{\FFF_2}$.

\begin{lemma}
The relation $\pi_k$ is a bijection for every $k\geq 1$.
\end{lemma}

\begin{proof}
Suppose that there are $z\not=y$ such that $xy\in\pi_k$ and $xz\in\pi_k$. Then, by~(\ref{OneLess}), there are two distinct sets $A,B\in\FFF_1$ such that
$$
k=r_{\FFF_1}(x)=|A|=|B|,\;\;A\setminus\{x\}\in\FFF_1,\;B\setminus\{x\}\in\FFF_1,
$$
and
$$
\aa(A)=\aa(A\setminus\{x\})+\{y\},\;\;\aa(B)=\aa(B\setminus\{x\})+\{z\}.
$$
We have
\begin{align*}
d\,(\aa(A),\aa(B))&=d\,(A,B)=d\,(A\setminus\{x\},B\setminus\{x\})\\
&=d\,(\aa(A)\setminus\{y\},\aa(B)\setminus\{z\}).
\end{align*}
Thus $y,z\in \aa(A)\cap\aa(B)$, that is, in particular, that $z\in \aa(A)\setminus\{y\}$, a contradiction, because $r_{\FFF_2}(z)=k$ and $|\aa(A)\setminus\{y\}|=k-1$.

By applying the above argument to $\aa^{-1}$, we prove that $\pi_k$ is a bijection.
\end{proof}

It follows from the previous lemma that $\pi$ is a permutation on $X$.

\begin{lemma}
$\alpha(A)=\pi(A)$ for any $A\in\FFF_1$.
\end{lemma}

\begin{proof}
We prove this statement by induction on $k=|A|$. The case $k=1$ is trivial, since $\aa(\{x\})=\{\pi_1(x)\}$ for $\{x\}\in\mathcal{F}_1$.

Suppose that $\aa(A)=\pi(A)$ for all $A\in\FFF_1$ such that $|A|<k$. Let $A$ be a set in $\FFF_1$ of cardinality $k$. By the wellgradedness property, there is a nested sequence $\{A_i\}_{0\leq i\leq k}$ of distinct sets in $\FFF_1$ with $A_0=\es$ and $A_k=A$. Thus, $A=A_{k-1}\cup\{x\}$ for some $x\notin A_{k-1}$. Clearly, $m=r_{\FFF_1}(x)\leq k$.

If $m=k$, then $\aa(A)=\aa(A_{k-1})\cup\{\pi(x)\}=\pi(A)$, by the definition of $\pi$ and the induction hypothesis.

Suppose now that $m<k$. There is a set $B\in\FFF_1$ containing $x$ such that $|B|=m$. By the wellgradedness property, there is a nested sequence $\{B_i\}_{0\leq i\leq m}$ of distinct sets in $\FFF_1$ with $B_0=\es$ and $B_m=B$. We have $x\notin B_i$ for $i<m$, since $m=r_{\FFF_1}(x)$. Therefore, $B=B_{m-1}\cup\{x\}$. Clearly,
$$
B_{m-1}\cap A\SB B\SB B_{m-1}\cup A.
$$
By~(\ref{so1 inclusion}), we have
$$
\aa(B)\SB \aa(B_{m-1})\cup \aa(A).
$$
Thus, by the induction hypothesis,
$$
\pi(B_{m-1})\cup\{\pi(x)\}=\pi(B)\SB \pi(B_{m-1})\cup \aa(A).
$$
Hence, $\pi(x)\in\aa(A)$. Since $\aa(A)=\pi(A_{k-1})\cup\{y\}$ for $y\notin\pi(A_{k-1})$ and $x\notin A_{k-1}$, we have $y=\pi(x)$, that is, $\aa(A)=\pi(A)$.
\end{proof}

In summary, we have the following theorem.

\begin{theorem} \label{WG homogeneity}
Any isometry between two partial subcubes of $\HHH(X)$ can be extended to an isometry of the cube $\HHH(X)$.
\end{theorem}

\begin{remark}
{\rm In the case of a finite set $X$ the previous theorem is a consequence of Theorem~19.1.2 in \cite{mD97}.}
\end{remark}

In the line of our arguments which led to the proof of Theorem~\ref{WG homogeneity} we used two kinds of isometries of $\HHH(X)$: isometries that map elements of $\HHH(X)$ to the empty set, and isometries defined by permutations on $X$. It is not difficult to show that these isometries generate the isometry group of $\HHH(X)$.

\begin{theorem} \label{so1 isometry group of cube}
The isometry group of $\HHH(X)$ is generated by permutations on the set $X$ and functions
$$
\aa_A:S\mapsto S\Delta A,\quad S\in\HHH(X).
$$
\end{theorem}

\begin{proof}
Clearly, $\aa_A$ is an isometry of $\HHH(X)$ and $\aa_A(A)=\es$. A permutation $\pi$ on $X$ defines an isometry $\hat{\pi}:\HHH(X)\rightarrow\HHH(X)$ by
$$
\hat{\pi}(S)=\{\pi(x): x\in S\}.
$$

Let $\aa:\HHH(X)\rightarrow\HHH(X)$ be an isometry of $\HHH(X)$ and let $A=\aa^{-1}(\es)$. Then the isometry $\aa_A\circ\aa^{-1}$ fixes $\es\in\HHH(X)$ and therefore defines a permutation $\pi:X\rightarrow X$ (singletons are on the distance $1$ from $\es$). Let $\bb=\hat{\pi}^{-1}\circ\aa_A\circ\aa^{-1}$. Since $\aa_A\circ\aa^{-1}$ fixes $\es$, we have $\aa_A\circ\aa^{-1}(\{x\})=\{\pi(x)\}$ for any $x\in X$. Hence, $\bb(\{x\})=\{x\}$ for all $x\in X$. For $S\in\HHH(X)$, we have
$$
|\bb(S)|=d(\bb(S),\es)=d(S,\es)=|S|,
$$
since $\bb(\es)=\es$. For any $x\in \bb(S)$, we have
$$
d(\{x\},S)=d(\{x\},\bb(S))=|\bb(S)|-1=|S|-1,
$$
which is possible only if $x\in S$. Thus $\bb(S)\SB S$. The same argument shows that $S\SB\bb(S)$. Thus $\bb$ is the identity mapping. It follows that $\aa=\hat{\pi}^{-1}\circ\aa_A$.
\end{proof}

\section{Linear Media}

The representation theorem (Theorem~\ref{FiniteRepresentationTheorem}) is a powerful tool for constructing media. We illustrate an application of this theorem by constructing a medium of linear orderings on a given finite or infinite countable set $Z$.

Let $Z=\{a_1,a_2,\ldots\}$ be a fixed (finite or infinite) enumeration of elements of $Z$. This enumeration defines a particular irreflexive linear ordering on $Z$ that we will denote by $L_0$.

\begin{definition}
A binary relation $R$ on $Z$ is said to be \emph{locally finite} if there is $n\in\mathbb{N}$ such that the restriction of $R$ to $\{a_{n+1},a_{n+2},\ldots\}$ coincides with the restriction of $L_0$ to the same set.
\end{definition}

Let $\LLL\OOO$ be the set of all locally finite irreflexive linear orders on the set $Z$. Note that if $Z$ is a finite set, then $\LLL\OOO$ is the set of all linear orderings on $Z$.

As usual, for a given $L\in\LLL\OOO$, we say that $x$ \emph{covers} $y$ in $L$ if $yx\in L$ and there is no $z\in Z$ such that $yz\in L,zx\in L$. Here and below $xy$ stands for an ordered pair of elements $x,y\in Z$. In what follows all binary relations are assumed to be locally finite for a given enumeration of $Z$.

\begin{lemma} \label{LO1}
Let $L$ be a linear order on $Z$. Then $L'=(L\setminus yx)\cup xy$ is a linear order if and only if $x$ covers $y$ in $L$.
\end{lemma}

\begin{proof}
Suppose $L'$ is a linear order and there is $z$ such that $yz\in L$ and $zx\in L$. Then $yz\in L'$ and $zx\in L'$ implying $yx\in L'$, a contradiction.

Suppose $x$ covers $y$ in $L$. Let $uv\in L'$ and $vw\in L'$. We need to show that $uw\in L'$. There are three possible cases.
\begin{enumerate}
	\item $uv=xy$ and $vw\in L,vw\not=yx$. Then $yw=vw\in L$ which implies $uw=xw\in L$, since $x$ covers $y$ in $L$. We have $uw\in L'$, since $uw=xw\not=yx$.
	\item $vw=xy$ and $uv\in L,uv\not=yx$. Then $ux=uv\in L$ which implies $uw=uy\in L$, since $x$ covers $y$ in $L$. We have $uw\in L'$, since $uw=uy\not=yx$.
	\item $uv\in L,uv\not=yx$ and $vw\in L,vw\not=yx$. Then $uw\in L$ and $uw\not= yx$, since $x$ covers $y$ in $L$. Therefore, $uw\in L'$.
\end{enumerate}
\end{proof}

We shall also need the following fact.

\begin{lemma} \label{LO2}
Let $P,Q$ and $R$ be complete asymmetric binary relations on $Z$. Then 
\begin{equation*}
P\cap R = Q\cap R\quad\eq\quad P=Q.
\end{equation*}
\end{lemma}

\begin{proof}
Suppose that $P\cap R=Q\cap R$ and let $xy\in P$. If $xy\in R$, then $xy\in Q$. Otherwise, $yx\in R$. Since $yx\notin P$, we have $yx\notin Q$ implying $xy\in Q$. Thus $P\SB Q$. By symmetry, $P=Q$.
\end{proof}

For $L\in\LLL\OOO$, we define
\begin{equation*}
\aa : L \mapsto L\cap L_0
\end{equation*}
By Lemma~\ref{LO2}, $\aa$ is a one--to--one mapping from $\LLL\OOO$ onto the set $\aa(\LLL\OOO)$ of partial orders.

Note that for any two locally finite binary relations $R$ and $Q$ on $Z$, the symmetric difference $R\Delta Q$ is a finite set. Thus the distance $d(R,Q)=|R\Delta Q|$ is a finite number. We use this fact in the proof of the following theorem.

\begin{theorem} \label{WGfamilyLO}
The family $\aa(\LLL\OOO)$ is a well graded family of subsets of $L_0$.
\end{theorem}
\begin{proof}
Let $P,P'$ be two distinct partial orders in $\aa(\LLL\OOO)$ and $L,L'$ be corresponding linear orders. It is easy to see that there is a pair $xy\in L$ such that $y$ covers $x$ and $xy\notin L'$. By Lemma~\ref{LO1}, $L''$ defined by
$$
L''=(L\setminus xy)\cup yx
$$
is a linear order. Then
$$
P''=L''\cap L_0 = [(L\cap L_0)\setminus(L_0\cap xy)]\cup(L_0\cap yx)=\begin{cases}
	P\setminus xy, & \text{if $xy\in L_0$,} \\
	P\cup yx, & \text{if $xy\notin L_0$,}
\end{cases}
$$
where $xy\in P$ if $xy\in L_0$ and $yx\notin P$ if $xy\notin L_0$. Hence, $P''\not= P$ and $d(P,P'')=1$. Clearly,
\begin{equation*}
L\cap L'\SB L''\SB L\cup L'.
\end{equation*}
Therefore
\begin{equation*}
P\cap P'=L\cap L'\cap L_0\SB P''=L''\cap L_0\SB (L\cup L')\cap L_0=P\cup P',
\end{equation*}
that is, $P''$ lies between $P$ and $P'$. Thus, 
\begin{equation*}
d(P,P')=d(P,P'')+d(P'',P')=1+d(P'',P')
\end{equation*}
and the result follows by induction.
\end{proof}

Since $\FFF=\aa(\LLL\OOO)$ is a well graded family of subsets of $X=L_0$, it is the set of states of the medium $(\FFF,\GGG_\FFF)$ with tokens defined by
\begin{equation*}
P\rho_{xy} = \begin{cases}
	P\cup xy, & \text{if $P\cup xy\in\aa(\LLL\OOO)$,} \\
	P, & \text{otherwise,}
\end{cases}
\end{equation*}
and
\begin{equation*}
P\tilde{\rho}_{xy} = \begin{cases}
	P\setminus xy, & \text{if $P\setminus xy\in\aa(\LLL\OOO)$,} \\
	P, & \text{otherwise,}
\end{cases}
\end{equation*}
for $xy\in L_0$ and $P=L\cap L_0\in\aa(\LLL\OOO)$.

We have
\begin{align*}
(L\cap L_0)\rho_{xy} & = \begin{cases}
	(L\cap L_0)\cup xy, & \text{if $(L\cap L_0)\cup xy\in\aa(\LLL\OOO)$,} \\
	L\cap L_0, & \text{otherwise,}
\end{cases} \\
& = \begin{cases}
	(L\cup xy)\cap L_0, & \text{if $(L\cup xy)\cap L_0=L'\cap L_0$,} \\
	L\cap L_0, & \text{otherwise,}
\end{cases}
\end{align*}
where $L'$ is some linear order. Since $yx\notin L_0$, we have 
\begin{equation*}
(L\cup xy)\cap L_0=[(L\setminus yx)\cup xy]\cap L_0 =L'\cap L_0.
\end{equation*}
By Lemma~\ref{LO2}, $(L\setminus yx)\cup xy$ is a linear order. We define
\begin{equation*}
L\tau_{xy}=\begin{cases}
	(L\setminus yx)\cup xy & \text{if $x$ covers $y$ in $L$,} \\
	L & \text{otherwise.}
\end{cases}
\end{equation*}
Then, for $xy\in L_0$,
\begin{equation*}
(L\cap L_0)\rho_{xy}=L\tau_{xy}\cap L_0.
\end{equation*}
A similar argument shows that, for $xy\in L_0$,
\begin{equation*}
(L\cap L_0)\tilde{\rho}_{xy}=L\tau_{yx}\cap L_0=L\tilde{\tau}_{xy}\cap L_0.
\end{equation*}

We obtained the set of tokens $\TTT=\{\tau_{xy}\}_{xy\in L_0}$ by `pulling back' tokens from the set $\GGG_\FFF$. The medium $(\LLL\OOO,\TTT)$ is isomorphic to the medium $(\FFF,\GGG_\FFF)$. In the case of a finite set $Z$, it is the \emph{linear medium} introduced in \cite{jF97}.

Simple examples show that $\aa(\LLL\OOO)$ is a proper subset of the set of all partial orders contained in $L_0$. In the case of a finite set $Z$ this subset is characterized in the following theorem.

\begin{theorem}
Let $L$ be a linear order on a finite set $Z$ and $P\SB L$ be a partial order. Then $P=L\cap L'$, where $L'$ is a linear order, if and only if $P'=L\setminus P$ is a partial order.
\end{theorem}

\begin{proof}
(1) Suppose $P=L\cap L'$. It suffices to prove that $P'=L\setminus P$ is transitive. Let $(x,y),\,(y,z)\in P'$. Then $(x,z)\in L$. Suppose $(x,z)\notin P'$. Then $(x,z)\in P$, implying $(x,z)\in L'$. Since $(x,y),\,(y,z)\in P'$, we have $(x,y)\notin L'$ and $(y,z)\notin L'$, implying $(y,x),\,(z,y)\in L'$, implying $(z,x)\in L'$, a contradiction.

(2) Suppose now that $P$ and $P'=L\setminus P$ are partial orders. We define $L'=P\cup {P'}^{-1}$ and prove that thus defined $L'$ is a linear order.

Clearly, relations $P,P^{-1},P',{P'}^{-1}$ form a partition of $(Z\times Z)\setminus\Delta$. It follows that $L'$ is a complete and antisymmetric binary relation.

To prove transitivity, suppose $(x,y),\,(y,z)\in L'$. It suffices to consider only two cases:

(i) $(x,y)\in P,\;(y,z)\in {P'}^{-1}$. Suppose $(x,z)\notin L'$. Then $(z,x)\in L'$. Suppose $(z,x)\in P$. Since $(x,y)\in P$, we have $(z,y)\in P$, a contradiction, since $(z,y)\in P'$. Suppose $(z,x)\in {P'}^{-1}$. Then $(x,z)\in P'$ and $(z,y)\in P'$ imply $(x,y)\in P'$, a contradiction. Hence, $(x,z)\in L'$.

(ii) $(y,z)\in P,\;(x,y)\in {P'}^{-1}$. Suppose $(x,z)\notin L'$. Then, again, $(z,x)\in L'$. Suppose $(z,x)\in P$. Since $(y,z)\in P$, we have $(y,x)\in P$, a contradiction, since $(y,x)\in P'$. Suppose $(z,x)\in {P'}^{-1}$. Then $(x,z)\in P'$ and $(y,x)\in P'$ imply $(y,z)\in P'$, a contradiction. Hence, $(x,z)\in L'$.

Clearly $P=L\cap L'$.
\end{proof}

We conclude this section with a geometric illustration of Theorem~\ref{WGfamilyLO}.

{\begin{figure}[h!]
\centerline{\includegraphics{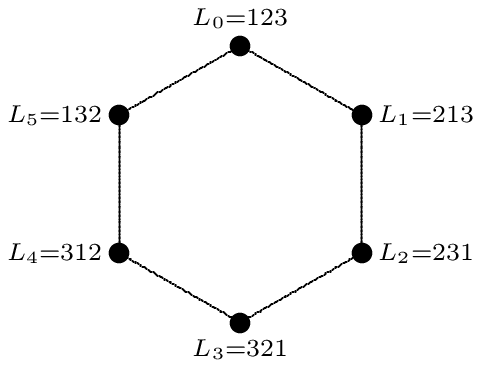}}
\caption{The diagram of $\LLL\OOO$.} \label{diagram LO}
\end{figure}
}

\begin{example}
{\rm Let $Z=\{1,2,3\}$. We represent linear orders on $Z$ by $3$--tuples. There are $3!=6$ different linear orders on $Z$:
\begin{equation*}
L_0=123,\quad L_1=213,\quad L_2=231,\quad L_3=321,\quad L_4=312,\quad L_5=132.
\end{equation*}
These relations are represented by the vertices of the diagram in Figure~\ref{diagram LO}.

One can compare this diagram with the diagram shown in Figure~5 in~\cite{jF97}.
The elements of $\FFF=\aa(\LLL\OOO)$ are subsets of $X=L_0$:
\begin{gather*}
L_0=\{12,13,23\},\quad L_1\cap L_0=\{13,23\},\quad L_2\cap L_0=\{23\}, \\
L_3\cap L_0=\es,\quad L_4\cap L_0=\{12\},\quad L_5\cap L_0=\{12,13\}.
\end{gather*}
These sets are represented as vertices of the cube on the set $L_0=\{12,13,23\}$ as shown in Figure~\ref{(LO,T)}.
}
\end{example}

{\begin{figure}[h!]
\centerline{\includegraphics{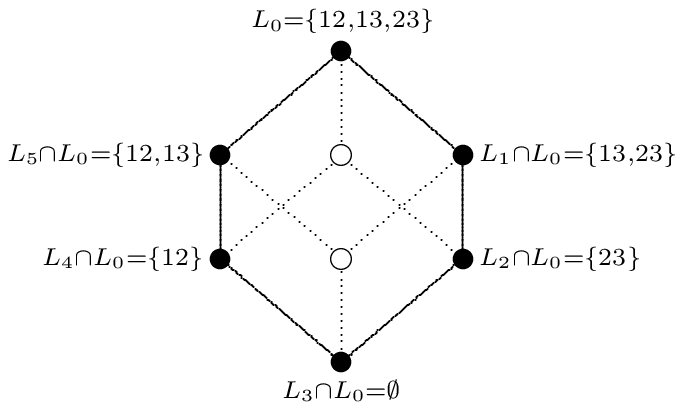}}
\caption{Partial cube representing $(\LLL\OOO,\TTT)$.} \label{(LO,T)}
\end{figure}
}

\section{Hyperplane arrangements}

In this section we consider an example of a medium suggested by Jean--Paul Doignon (see Example~2 in~\cite{jFsO02}).

Let $\mathcal{A}$ be a locally finite arrangements of affine hyperplanes in $\mathbb{R}^r$, that is a family of hyperplanes such that any open ball in $\mathbb{R}^r$ intersects only finite number of hyperplanes in $\mathcal{A}$~\cite[Ch.~V,\;\S 1]{nB02}. Clearly, there are only countably many hyperplanes in $\mathcal{A}$, so we can enumerate them, $\mathcal{A}=\{H_1,H_2,\ldots\}$. Every hyperplane is given by an affine linear function $\ell_i(\boldsymbol{x})=\sum_{j=1}^r a_{ij}x_j+b_i$, that is, $H_i=\{\boldsymbol{x}\in\mathbb{R}^r : \ell_i(\boldsymbol{x})=0\}$.

In what follows, we construct a token system $(\SSS,\TTT)$ associated with an arrangement $\mathcal{A}$ and show that this system is a medium.

We define the set $\SSS$ of states to be the set of connected components of $\mathbb{R}^r\setminus\cup\,\mathcal{A}$. These components are called \emph{regions}~\cite{aB99} or \emph{chambers}~\cite{nB02} of $\mathcal{A}$. Each state $P\in\SSS$ is an interior of an $r$--dimensional polyhedron in $\mathbb{R}^r$.

To every hyperplane in $\mathcal{A}$ corresponds an ordered pair $(H,H')$ of open half spaces $H$ and $H'$ separated by this hyperplane. This ordered pair generates a transformation $\tau_{H,H'}$ of the states. Applying $\tau_{H,H'}$ to some state $P$ results in some other state $P'$ if $P\SB H,\;P'\SB H'$ and regions $P$ and $P'$ share a facet which is included in the hyperplane separating $H$ and $H'$; otherwise, the application of $\tau_{H,H'}$ to $P$ does not change $P$. We define the set $\TTT$ of tokens to be the set of all $\tau_{H,H'}$. Clearly, $\tau_{H,H'}$ and $\tau_{H',H}$ are reverses of each other.

\begin{theorem} \label{HyperplaneMedia}
$(\SSS,\TTT)$ is a medium.
\end{theorem}

\begin{proof}
In order to prove that $(\SSS,\TTT)$ is a medium, we show that it is isomorphic to a medium of a well graded family of sets. Let $J=\{1,2,\ldots\}$. We denote $H_i^+=\{\boldsymbol{x}\in\mathbb{R}^r : \ell_i(\boldsymbol{x})>0\}$ and $H_i^-=\{\boldsymbol{x}\in\mathbb{R}^r : \ell_i(\boldsymbol{x})<0\}$, open half spaces separated by $H_i$. Each region $P$ is an intersection of open half spaces corresponding to hyperplanes in $\mathcal{A}$. We define $J_P = \{j\in J : P\SB H_j^+\}$. Clearly, $P\mapsto J_P$ defines a bijection from $\SSS$ to $\SSS'=\{J_P:P\in\SSS\}$. It is also easy to see that $\cap\,\SSS'=\es$ and $\cup\,\SSS'=J$.

Given $k\in J$, we define transformations $\tau_k$ and $\tilde{\tau}_k$ of $\SSS'$ as follows:
$$
	J_P\tau_k = \begin{cases}
		J_P\cup\{k\} &\text{if $H_k$ defines a facet of $P$,} \\
		J_P &\text{otherwise},
\end{cases}
$$
and
$$
	J_P\tilde{\tau}_k = \begin{cases}
		J_P\setminus\{k\} &\text{if $H_k$ defines a facet of $P$}, \\
		J_P &\text{otherwise}.
\end{cases}
$$
Let $P$ be a region of $\mathcal{A}$ and let $H_k$ be a hyperplane in $\mathcal{A}$ defining a facet of $P$. There is a unique region $P'$ sharing this facet with $P$. Moreover, $H_k$ is the only hyperplane separating $P$ and $P'$. It follows that $J_P\tau_k=J_{P'}$ if $k\notin J_P$ and $J_P\tilde{\tau}_k=J_{P'}$ if $k\in J_P$. Thus transformations $\tau_k$ and $\tilde{\tau}_k$ are well defined.

We denote $\TTT'$ the set of all transformations $\tau_i,\;\tilde{\tau}_i,\;i=1,\ldots,n$. Clearly, the correspondences $\tau_{H_i^+,H_i^-}\mapsto\tilde{\tau}_i$ and $\tau_{H_i^-,H_i^+}\mapsto\tau_i$ define a bijection from $\TTT$ to $\TTT'$. This bijection together with the bijection from $\SSS$ to $\SSS'$ given by $P\mapsto J_P$ define an isomorphism of two token systems, $(\SSS,\TTT)$ and $(\SSS',\TTT')$.

It remains to show that $\SSS'$ is a well graded family of subsets of $J$.

Clearly, $k\in J_P\Delta J_Q$ if and only if $H_k$ separates $P$ and $Q$. Since $\mathcal{A}$ is locally finite, there is a finite number of hyperplanes in $\mathcal{A}$ that separate two regions. Thus $J_P\Delta J_Q$ is a finite set for any two regions $P$ and $Q$. Let $d$ be the usual Hamming distance on $\SSS'$, i.e, $d(J_P,J_Q)=|J_P\Delta J_Q|$. Thus $d(J_P,J_Q)$ is equal to the number of hyperplanes in $\mathcal{A}$ separating $P$ and $Q$.

Let $\boldsymbol{p}\in P$ and $\boldsymbol{q}\in Q$ be points in two distinct regions $P$ and $Q$. The interval $[\boldsymbol{p},\boldsymbol{q}]$ has a single intersection point with any hyperplane separating $P$ and $Q$. Moreover, a simple topological argument shows that we can always choose $\boldsymbol{p}$ and $\boldsymbol{q}$ in such a way that different hyperplanes separating $P$ and $Q$ intersect $[\boldsymbol{p},\boldsymbol{q}]$ in different points. Let us number these points in the direction from $\boldsymbol{p}$ to $\boldsymbol{q}$ as follows
\begin{equation*}
\boldsymbol{r}_0 = \boldsymbol{p}, \boldsymbol{r}_1,\ldots,\boldsymbol{r}_{k+1}=\boldsymbol{q}.
\end{equation*}
Each open interval $(\boldsymbol{r}_i,\boldsymbol{r}_{i+1})$ is an intersection of $[\boldsymbol{p},\boldsymbol{q}]$ with some region which we denote $R_i$ (in particular, $R_0=P$ and $R_k=Q$). Moreover, by means of this construction, points $\boldsymbol{r}_i$ and $\boldsymbol{r}_{i+1}$ belong to facets of $R_i$. We conclude that regions $R_i$ and $R_{i+1}$ are adjacent, that is, share a facet, for all $i=0,\ldots,k-1$. Clearly, $d(J_{P_i},J_{P_{i+1}})=1$ for all $i=0,\ldots,k-1$ and $d(J_P,J_Q)=k$. Thus, $\SSS'$ is a well graded family of subsets of $J$.
\end{proof}

The \emph{region graph} $G$~\cite{aB99} of the arrangement $\mathcal{A}$ has $\SSS$ as the set of vertices; edges of $G$ are pairs of adjacent regions in $\SSS$. It follows from Theorem~\ref{HyperplaneMedia} that $G$ is a partial cube.

In the case of a finite arrangement $\mathcal{A}$, the graph $G$ is the \emph{tope graph} of the oriented matroid associated with the arrangement $\mathcal{A}$. It follows from Proposition~4.2.3 in~\cite{aB99} that $G$ is an isometric subgraph  of the $n$--cube, where $n$ is the number of hyperplanes in $\mathcal{A}$. Thus our Theorem~\ref{HyperplaneMedia} is an infinite dimensional analog of this result.

To give geometric examples of infinite partial cubes, let us consider locally finite line arrangements $\mathcal{A}$ in the plane $\mathbb{R}^2$. The closures of the regions of a given $\mathcal{A}$ form a tiling~\cite{bG87,mS95} of the plane. The region graph of this tiling is the $1$--skeleton of the dual tiling.

\begin{example}
{\rm Let us consider a line arrangement $\mathcal{A}$ shown in Figure~\ref{hex mosaic} by dotted lines.
The regions of this line arrangement are equilateral triangles that form $(3^6)$ mosaic (an edge--to--edge planar tiling by regular polygons; for notations and terminology see, for instance,~\cite{mD02,bG87}). The $1$--skeleton of the orthogonally dual~\cite{mS95} mosaic $(6^3)$ is the region graph of $\mathcal{A}$. This graph is also known as the hexagonal lattice in the plane.
By Theorem~\ref{HyperplaneMedia}, the hexagonal lattice is an infinite partial cube. This lattice is isometrically embeddable into the graph of the cubical lattice $\Zee^3$~\cite{mD02}.
}
\end{example}

{\begin{figure}[h!]
\centerline{\includegraphics{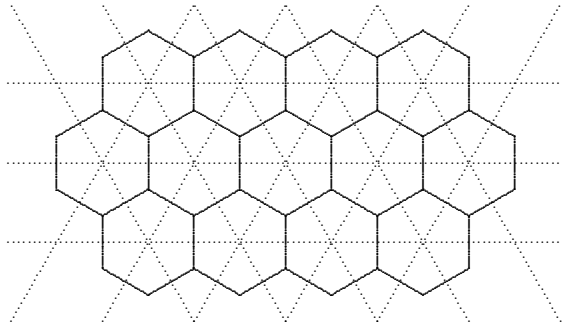}}
\caption{Hexagonal lattice ($(6^3)$ mosaic).} \label{hex mosaic}
\end{figure}
}

\begin{example}
{\rm 
Another example of an infinite partial cube is shown in Figure~\ref{truncated mosaic}. There, the region graph is the $1$--skeleton of $(4.8^2)$ mosaic also known~\cite{mD02} as the truncated net $(4^4)$. Like in the previous case, this mosaic is orthogonally dual to the tiling defined by the line arrangement shown in Figure~\ref{truncated mosaic} and the region graph can be isometrically embedded into $\Zee^4$~\cite{mD02}.
}
\end{example}

{\begin{figure}[h!]
\centerline{\includegraphics{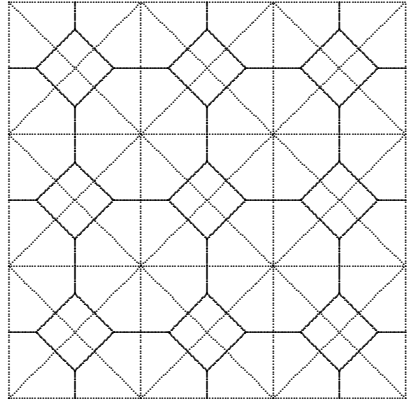}}
\caption{$(4.8)$ mosaic.} \label{truncated mosaic}
\end{figure}
}

\begin{example}
{\rm 
A more sophisticated example of an infinite partial cube was suggested by a referee. This is one of the Penrose rhombic tilings (see, for instance,~\cite{dB81,mS95}) a fragment of which is shown in Figure~\ref{Penrose}~\cite[Ch.~9]{sW99}. The construction suggested by de~Bruijn~\cite{dB81} demonstrates that the graph of this tiling is the region graph of a particular line arrangement known as a pentagrid. This graph is isometrically embeddable in $\Zee^5$~\cite{dB81,mD02}.
}
\end{example}

{\begin{figure}[h!]
\centerline{\includegraphics{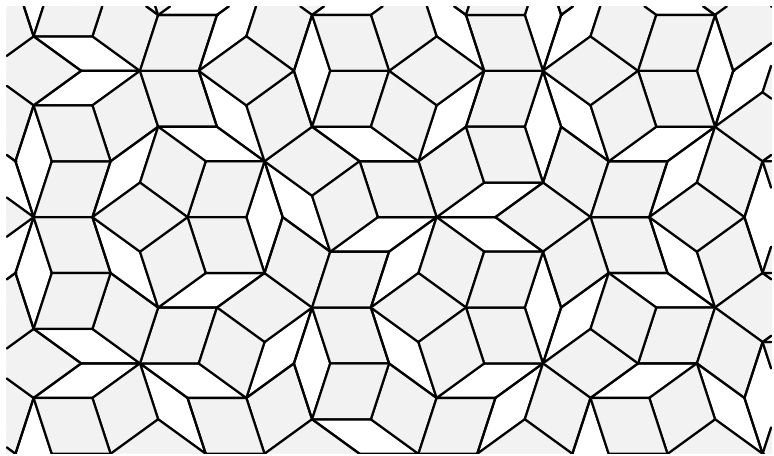}}
\caption{A Penrose rhombic tiling.} \label{Penrose}
\end{figure}
}

\section*{Acknowledgments}
The author is grateful to Jean--Claude Falmagne for his careful reading of the original manuscript and many helpful suggestions, and to Jean--Paul Doignon for his comments on the results presented in Section~7. I also thank the referees for their constructive criticism.

\end{document}